\documentclass[leqno,11pt]{article}
\usepackage{amssymb,latexsym,mathrsfs}
\usepackage[dvips]{graphicx}
\usepackage{url}

 \catcode`@=11\def\@oddhead{\hbox{}\hfil\rm\thepage}\def\@oddfoot{}
 \def\@evenhead{\hbox{}\hfil\rm\thepage}\def\@evenfoot{}
 \@twosidetrue \catcode`@=12
\newtheorem{prp}{Proposition}
\newtheorem{lem}[prp]{Lemma}\newtheorem{thm}[prp]{Theorem}
\newtheorem{cor}[prp]{Corollary}
\newenvironment{prf}{\begin{trivlist}\item[\emph{Proof.}]}{\end{trivlist}
  \medskip\par}
\newenvironment{prfof}[1]{\begin{trivlist}\item[\emph{Proof of #1.}]}{
  \end{trivlist} \medskip \par}
\newenvironment{rem}{\begin{trivlist}\item[\emph{Remark.}]}{\end{trivlist}
  \medskip\par}
\def\prpb{\begin{prp}}\def\prpe{\end{prp}}
\def\lemb{\begin{lem}}\def\leme{\end{lem}}
\def\thmb{\begin{thm}}\def\thme{\end{thm}}
\def\corb{\begin{cor}}\def\core{\end{cor}}
\def\prfb{\begin{prf}}\def\prfe{\end{prf}}
\def\prfofb#1{\begin{prfof}{#1}}\def\prfofe{\end{prfof}}
\def\remb{\begin{rem}}\def\reme{\end{rem}}
\def\prpa#1{\label{p:#1}}\def\prpu#1{Proposition~\ref{p:#1}}
\def\lema#1{\label{l:#1}}\def\lemu#1{Lemma~\ref{l:#1}}
\def\thma#1{\label{t:#1}}\def\thmu#1{Theorem~\ref{t:#1}}
\def\cora#1{\label{c:#1}}
\def\seca#1{\label{s:#1}}\def\secu#1{\S~\ref{s:#1}}

\def\itmb{\begin{enumerate}}\def\itme{\end{enumerate}}
\def\itdb{\begin{itemize}}\def\itde{\end{itemize}}
\def\ittb{\begin{description}}\def\itte{\end{description}}
\def\arrb#1{\begin{array}{#1}}\def\arre{\end{array}}
\def\tabb#1{\par\noindent\begin{tabular}{#1}}
\def\tabe{\end{tabular}\par\noindent}
\def\eqna#1{\label{e:#1}}\def\eqnu#1{(\ref{e:#1})}

\def\QED{\relax\ifmmode\let\@tempa\relax\ifcase\@eqcnt\def\@tempa{& & &}\or
  \def\@tempa{& &}\else\def\@tempa{&}\fi\@tempa $\Box$ \else\hfill $\Box$ \fi}
\def\DDD{\relax\ifmmode\let\@tempa\relax\ifcase\@eqcnt\def\@tempa{& & &}\or
 \def\@tempa{& &}\else\def\@tempa{&}\fi\@tempa $\Diamond$
 \else\hfill $\Diamond$ \fi}

\def\Rom#1{\uppercase\expandafter{\romannumeral#1}}
\def\dsp{\displaystyle}
\def\eps{\epsilon}
\def\AHFA{\hbox{\vrule height16pt depth10pt width0pt}\displaystyle}
\def\limf#1{\displaystyle \lim_{#1\to\infty}}

\def\limsup{\displaystyle \mathop{\overline{\lim}}\limits}

\def\norm#1{\dsp\left\| #1 \right\|}

\def\Ccomb#1#2{\setbox0=\hbox{$\displaystyle\mathrm{C}$}\setbox1=\hbox{%
$\scriptstyle #1$}\kern \wd1{\mathrm{C}}_{\kern -1.05\wd0\kern -0.99\wd1{#1}
 \kern 1.15\wd0{#2}}}

\def\vec#1{{\mathbf{#1}}} \def\vec#1{\overrightarrow{#1}}
\def\clvec#1#2#3{\def\clvecone{#3}\left(\arrb{c} \dsp #1\\ \dsp #2
 \ifx\clvecone\empty\else\\ \dsp #3\fi\arre\right)}

\def\pderiv#1#2{\dsp\frac{\partial\,#1}{\partial#2}}

\def\bar#1{\overline{#1}}
\def\le{\leqq} \def\ge{\geqq} \def\geq{\geqq}

\def\reals{{\mathbb R}}
\def\nreals#1{{\mathbb R}^{#1}}
\def\preals{[0,\infty)} 
\def\pintegers{{\mathbb Z}_+}
\def\nintegers{{\mathbb N}}

\def\prb#1{\def\prbone{#1}
  \ifx\prbone\empty{\mathrm{P}}\else{\mathrm{P[\;}}#1{\mathrm{\;]}}\fi}
\def\prbseq#1#2{\def\prbseqone{#2}
  \ifx\prbseqone\empty{\mathrm{P}}_{#1}\ignorespaces
  \else{\mathrm{P}}_{#1}{\mathrm{[\;}}#2{\mathrm{\;]}}\fi}
\def\EE#1{{\mathrm{E[\;}}#1{\mathrm{\;]}}}
\def\EEseq#1#2{\def\EEseqone{#2}
  \ifx\EEseqone\empty{\mathrm{E}}_{#1}\else
 {\mathrm{E}}_{#1}{\dsp\mathrm{[\;}}#2{\mathrm{\;]}}\fi}

\def\VVseq#1#2{\def\VVseqone{#2}
  \ifx\VVseqone\empty{\matrm{V}}_{#1}\else
 {\mathrm{V}}_{#1}{\dsp\mathrm{[\;}}#2{\mathrm{\;]}}\fi}

\def\Ccomb#1#2{\setbox0=\hbox{$\displaystyle\mathrm{C}$}\setbox1=\hbox{%
$\scriptstyle #1$}\kern \wd1{\mathrm{C}}_{\kern -1.05\wd0\kern -0.99\wd1{#1}
 \kern 1.15\wd0{#2}}}

\def\ssN{^{(N)}}

\def\ssNtheta{^{(N,\theta)}}
\def\ssNyC{^{(N,y_C)}}

\def\wnorm#1{\dsp\left\| #1 \right\|_{\rm{T}}}
\def\chrfcn#1{\mathop{\mathbf{1}}\nolimits_{#1}}

\title{
Cancellation of fluctuation in stochastic ranking process with space-time dependent intensities,
}
\author{
Tetsuya Hattori
\\ 
\small Laboratory of Mathematics, Faculty of Economics, Keio University, 
\\
\small Hiyoshi Campus, 4--1--1 Hiyoshi, Yokohama 223-8521, Japan
\\ \small URL: \url{http://web.econ.keio.ac.jp/staff/hattori/research.htm}
\\ \small email: \url{hattori@econ.keio.ac.jp}
} 
\date{\today}

\begin{document}

\maketitle

\begin{abstract} 
We consider the stochastic ranking process with space-time dependent unbounded jump rates for the particles. We prove that the joint empirical distribution of jump rate and scaled position converges almost surely to a deterministic distribution in the infinite particle limit. We assume topology of weak convergence for the space of distributions, which implies that the fluctuations among particles with different jump rates cancel in the limit. The results are proved by first finding an auxiliary stochastic ranking process, for which a strong law of large numbers is applied, and then applying a multi time recursive Gronwall's inequality. The limit has a representation in terms of non-Markovian processes which we call point processes with last-arrival-time dependent intensities.
\end{abstract}

\noindent 
2000 \textit{Mathematics Subject Classification.}
Primary 60K35;
Secondary 
82C22.

\noindent 
\textit{Key words}. 
Stochastic ranking process, hydrodynamic limit,
law of large numbers, complete convergence, 
last-arrival-time dependent intensity, Gronwall inequality.

\newpage

\setcounter{section}{0}

\setcounter{equation}{0}

\section{Introduction.}
\seca{1}

Let $N$ be a positive integer and $T>0$.
$T$ is an arbitrary constant fixed throughout the paper,
and we are interested in the limit $N\to\infty$.
A stochastic ranking process is a stochastic system of $N$ particles
on a line segment $[0,1]$ for a time interval $[0,T]$,
defined as follows.
Let $W\subset C^1([0,1]\times[0,T];\preals)$ be a set of 
non-negative valued $C^1$ functions 
in two variables such that the partial derivative with respect to the
first variable is bounded on $[0,1]\times[0,T]$.
We write
\begin{equation}
\eqna{Wnorm}
\wnorm{f}=\sup_{(z,s)\in[0,1]\times[0,T]} |f(z,s)|,
\end{equation}
for a function $f:\ [0,1]\times[0,T]\to \reals$, and put
\begin{equation}
\eqna{Cw}
C_W= \sup_{w\in W} \wnorm{\pderiv{w}{y}} <\infty.
\end{equation}
Let $w_1$, $w_2$, $\ldots$ be an infinite sequence in $W$.

Let $y\ssN_1$, $y\ssN_2$, $\ldots$, $y\ssN_N$ be a permutation of
$\dsp\{\frac{i}N\mid i=0,1,\ldots,N-1\}$.
Then a stochastic ranking process is 
a system of stochastic processes $\{Y\ssN_i\mid i=1,2,\ldots,N\}$
defined on a probability space $(\Omega,\mathcal{F},\prb{})$ by
\begin{equation}
\eqna{SRPY}
\begin{array}{l}\dsp
 Y\ssN_i(t) = y\ssN_i 
 \\ \dsp \phantom{Y=}
+ \frac1N\sum _{j=1} ^N
 \int _{s\in (0,t]} \int _{\xi \in \preals }
 \chrfcn{ Y\ssN_j(s-)>Y\ssN_i(s-) }
 \chrfcn{\xi\in [0,w_{j} (Y\ssN_j(s-),s))}
 \nu\ssN _j(d\xi ds)
 \\ \dsp \phantom{Y=}
 - \int _{s\in (0,t]}  \int _{\xi \in \preals } Y\ssN_i(s-)
 \chrfcn{\xi\in [0,w_i(Y\ssN_i(s-),s))} 
 \nu\ssN _i(d\xi ds), 
\\ \dsp i=1,2,\dots ,N,\ t\geq 0.
 \arre
\end{equation}
Here, $\chrfcn{A}$ is the indicator function of an event $A$,
and for each $N$, $\nu\ssN_i$, $i=1,2,\ldots,N$, are independent
Poisson random measures on $\preals \times\preals $
with uniform unit intensity measures, i.e.,
$\EE{\nu\ssN_i([a,b]\times[c,d])}=(b-a)(d-c)$ for $b>a\ge0$ and $d>c\ge0$,
and $\nu\ssN_i(A)$ and $\nu\ssN_i(B)$ are independent Poisson variables
if $A\cap B=\emptyset$.

If we put
\begin{equation}
\eqna{PPLATDPY}
\tilde{\nu}\ssN_i(t)=\int_{s\in(0,t]}\int _{\xi \in \preals }
 \chrfcn{\xi\in [0,w_{i} (Y\ssN_i(s-),s))} \nu\ssN _i(d\xi ds),
\end{equation}
then \eqnu{SRPY}, the position of the particle $i$ at time $t$, 
is expressed as
\[
 Y\ssN_i(t) = y\ssN_i 
+ \frac1N\sum _{j=1} ^N \int_0^t
 \chrfcn{ Y\ssN_j(s-)>Y\ssN_i(s-) } \tilde{\nu}\ssN _j(ds)
 - \int_0^t Y\ssN_i(s-) \tilde{\nu}\ssN _i(ds).
\]
We then see that the time development of $\{Y\ssN_i\}$
is determined by the move-to-front rules \cite{Tsetlin1963} driven by
the point processes $\{\tilde{\nu}\ssN_i\}$ with space-time
dependent intensities given by the `density functions' $\{w_i\}$.
If $w_i$ is constant then the process $\tilde{\nu}\ssN_i$
is the Poisson process, while in general, 
its increments depend on past and is dependent on other particles.
The move-to-front rule in particular implies that the particle
system $\{Y\ssN_i\}$ as a whole 
takes values in the rearrangement of $\dsp\frac{i}N$, $i=0,1,\ldots,N-1$.
Each particle either increases
its position by $\dsp \frac1N$, or else takes the value $0$, i.e.
jumps to the top position, as $t$ increases.

A starting point for our study is the joint empirical distribution of $w_i$
and the position, given by
\begin{equation}
\eqna{mut}
\mu\ssN_{t} =\frac1N \sum_{i=1}^N \delta_{(w_i,Y\ssN_i(t))}\,.
\end{equation}
Here $\delta_c$ is a unit measure concentrated at $c$.
(We will use this notation for a unit measure on any probability space.)
$\mu\ssN_{\cdot}$ is a stochastic process
taking values in the set of Borel probability measures,
with initial distribution being
\begin{equation}
\eqna{initialdistributionfcn}
\mu\ssN_0
=\frac1N \sum_{i=1}^N \delta_{(w_i,y\ssN_i)}.
\end{equation}

By considering a process $X\ssN_i(t)=N Y\ssN_i(t)+1$ taking
values in positive integers, we see that
a stochastic ranking process is a model of ranking system,
such as the sales ranks found at online bookstores\cite{HH071,HH072,timedep,
mv2frnt,ranking,dojin,Nagahata10,Nagahata11,posdep}.
As a model of popularity ranks of an online bookstore,
the move-to-front rule defines the rank as a stochastic number with
the `latest purchased book as most popular' rule.
That the intensity $w_i$ differs for different book $i$ represents that
there are popular books and less popular ones. (Many books on mathematics
perhaps provide examples of the latter.)
See also \cite{posdep} and references therein for more background.

We will prove existence of hydrodynamic limit,
that assuming convergence of $\mu\ssN_0$ as $N\to\infty$
we have convergence of $\mu\ssN_t$ for all $t\in[0,T]$.
The common standard quantities effective for the move-to-front rules
are the characteristic curves $Y\ssN_C$ defined by
\begin{equation}
\eqna{YNC2}
\arrb{l}\dsp 
Y\ssN_C(\gamma,t) = y_0
\\ \AHFA \dsp  \phantom{Y=}
 + \frac1N \sum_{j=1}^N \int _{s\in (t_0,t]} \int _{\xi \in \preals } 
\chrfcn{Y\ssN_j(s-)\ge Y\ssN_C(\gamma,s-)} 
\\ \dsp \phantom{Y= + \frac1N \sum_{j=1}^N \int }\times
\chrfcn{\xi \in [0,w_j(Y\ssN_j(s-),s))}\nu\ssN _j (d\xi ds),
\\ \dsp
\gamma=(y_0,t_0)\in [0,1]\times[0,T],\ t\ge t_0,
\arre
\end{equation}
and a set of spatial distribution functions $\varphi\ssN$
of $\mu\ssN_t$ defined by
\begin{equation}
\eqna{varphiN}
\arrb{l}\dsp
\varphi\ssN(dw,\gamma,t)=\mu\ssN_t(dw\times [Y\ssN_C(\gamma,t),1]),
\\ \dsp
\gamma=(y_0,t_0)\in [0,1]\times[0,T],\ t\ge t_0,\ w\in W.
\arre
\end{equation}
The latter is a refinement of the former in the sense that
\eqnu{mut} and \eqnu{varphiN} imply
\begin{equation}
\eqna{varphiNYNC}
\arrb{l}\dsp
Y\ssN_C(\gamma,t)=y_0+\frac{[N(1-y_0)]}{N} -\varphi\ssN(W,\gamma,t),
\\ \dsp 
\gamma=(y_0,t_0)\in[0,1]\times[0,T],\ t\ge t_0\,.
\arre
\end{equation}
Since $\varphi\ssN$ determines $\mu\ssN_t$,
the convergence problem of $\mu\ssN_t$ is reduced to that of $\varphi\ssN$.

If the intensity densities $w\in W$ are \textit{independent} of $y$,
then $\varphi\ssN$ is an arithmetic mean of \textit{independent} stochastic
processes, so that a strong law of large numbers for a sum of 
\textit{independent processes} can be applied to prove existence of 
almost sure $N\to\infty$ limit under reasonable assumptions \cite{timedep}.
In contrast, when the intensity densities depend on $y$,
the problem is a more involved one of law of large numbers for 
\textit{dependent processes}.
The case of spatially varying intensity densities  was first studied 
in \cite{posdep},  where we found the existence of the limit
under restricting assumptions.
One assumption was boundedness of $W$,
which excludes, for example, Pareto (power law) distributions
which may be of interest in applying stochastic ranking processes to 
social studies\cite{dojin,ranking}.
Another, a mathematically more essential assumption in \cite{posdep}
was that we adopted
total variation norm for the topology of the space of Borel measures.
This means that, among other points, 
the law of large numbers proved in \cite{posdep} is a cancellation of 
fluctuations among processes $\{Y\ssN_i\}$ having the same 
associated intensity densities $w_i=w$.
This assumption is to be compared with the results in \cite{timedep},
where we proved a corresponding convergence theorem (for the easier case
of spatially constant intensity densities) with topology 
of the space of Borel measures induced by weak convergence,
which, for example, allows all the $w_i$'s to be different
and the cancellation of fluctuation is still implied.
It would be mathematically interesting to see how this
law of large numbers (fluctuation cancellation) mechanism,
which escaped from our hands in \cite{posdep},
would be stable against introduction of the spatial dependence of $w\in W$.

The essential ingredient of the proof of this paper may be summarized 
as follows:
\begin{itemize}
\item
Discovery of the  point process with last-arrival-time dependent 
intensity \cite{fluid14,fluid14k}.
The hydrodynamic limit turned out to have an expression in terms of
probabilities of the process.
This expression was absent in \cite{posdep}, which was a partial cause of
a technical boundedness assumptions on $W$.
The process lack independent increment properties,
which is a remnant of stochastic dependence among the particles
through the position dependence of $w\in W$.
This is to be compared with the earlier studies for
position independent $w$, where 
the corresponding quantities have explicit formulas using
exponentials of integration of $w$, related to the probabilities of
the Poisson process \cite{HH071,timedep}.

\item
Discovery of an intermediate model, defined in \secu{fSRP},
which we call the \textit{flow driven stochastic ranking process}.
The `intensity densities' for the flow driven stochastic
ranking process are involved but without position dependence,
hence the distribution function $\varphi$ for this model
is a sum of independent processes,
and a standard law of large numbers
has a chance of explaining the fluctuation cancellation mechanism.
(Introduction of an intermediate model resembles a notion of 
local equilibrium which appears in the hydrodynamic limit 
for diffusions. The limit of the stochastic ranking process is,
in terms of fluid dynamics, a one-sided flow with evaporation
from upper stream $y=0$ to the down stream $y=1$, 
and the correspondence to diffusion is only a kind of metaphor.)

Since the stochastic dependence among the particles 
induced by the move-to-front rule is handled by introduction
of distribution functions $\varphi$, the real challenge from 
a viewpoint of mathematical analysis is the stochastic dependence 
through position dependence of $w\in W$,
which was first studied in \cite{posdep}.
In the reference we adopted a sophisticated method 
(than adopted in this paper) based on submartingale inequalities,
which worked well with the restricting assumptions in the reference, and
enabled a relatively quick proof without introducing intermediate models.

\item
Development and application of 
a uniform strong complete law of large numbers for 
independent monotone function valued random variables.

Since the flow driven stochastic ranking process
is stochastically similar to the stochastic ranking process with
position independent intensities,
a strong law of large numbers for independent processes
is applicable.
Since, however, we later need to compare this intermediate model with 
the original model,
we apply uniform convergence results stronger than in 
earlier works \cite{timedep},
where the law of large numbers for the independent processes
was practically the final goal.

\item
Application of a hierarchy of multi time Gronwall type inequalities.
The last step is to prove that the original 
stochastic ranking process has a same limit with the
(appropriately chosen) flow driven stochastic ranking process.
To evaluate the difference, we couple the two models,
and resort to a multi time variable recursive version of 
Gronwall type inequality which we develop in \secu{mGronwall}.

\end{itemize}

In \cite{posdep}, we also proved the occurrence of propagation of chaos.
Namely, for each integer $L$, the tagged particle system 
\[
(Y\ssN_1(t),Y\ssN_2(t),\ldots,Y\ssN_L(t))
\]
converges to a limit process uniformly in $t\in[0,T]$ as $N\to\infty$,
if the system of the initial positions 
$(y\ssN_1,y\ssN_2,\ldots,y\ssN_L)$ converges,
and the components of the limit are independent of each other.
A corresponding result is also proved in this paper.

The plan of the present paper is as follows.
In \secu{statement} we state the main results, that
the stochastic ranking process with space-time dependent intensities
has a limit characterized by 
the point processes with last-arrival-time dependent intensities.
For convenience, 
we will summarize the definition and relevant results of 
the point processes in \secu{ptproc}.
In \secu{fSRP} we formulate the flow driven stochastic ranking process
and prove the existence of the strong uniform law of large numbers
(the large particle number limit).
In \secu{mGronwall} we formulate and prove a hierarchy of a 
multi time version of the 
Gronwall's inequality, which we use in \secu{MainProof} to complete a proof 
of the main theorem stated in \secu{statement}.

\section{Formulation and the main results.}
\seca{statement}


To state the main results precisely,
we first formulate the quantities which appear in the
infinite particle limit of the stochastic ranking process.
Denote the sets of initial ($t=0$) points
in the space-time $[0,1]\times[0,T]$,
the set of upper stream boundary ($y=0$) points,
and their union, the set of initial/boundary points, respectively by
\begin{equation}
\eqna{Gamma}
\begin{array}{l}\dsp
\Gamma_b=\{(0,s)\mid 0\le s\le T\},
\\
\Gamma_i=\{(z,0)\mid 0\le z\le 1\},
\\
\Gamma=\Gamma_b\cup\Gamma_i\,.
\end{array}
\end{equation}
For $t\in[0,T]$, denote 
the set of initial/boundary points up to time $t$ by
\begin{equation}
\eqna{Gammat}
\Gamma_t=\{(z,t_0)\in \Gamma \mid t_0\le t\}
=\Gamma_i\cup \{(0,t_0)\in\Gamma_b \mid 0\le t_0\le t\},
\end{equation}
and the set of admissible pairs of the initial/boundary point $\gamma$
and time $t$ by
\begin{equation}
\eqna{domyC}
\Delta_T=\{(\gamma,t)\in\Gamma_T\times [0,T]\mid \gamma\in\Gamma_t\}.
\end{equation}

Define the set of \textit{flows} $\Theta_T$ on $[0,1]\times[0,T]$ by
\begin{equation}
\eqna{class_of_conti_flows}
\begin{array}{l}\dsp
\Theta_T=\{\theta:\ \Delta_T\to[0,1]\mid \mbox{continuous, }
\theta((y_0,t_0),t_0)=y_0,\ (y_0,t_0)\in\Gamma_T,
\\ \dsp \phantom{\Theta_T=\{} \mbox{ %
surjective and non-increasing in $\gamma$ for each $t$,
 }\\ \dsp \phantom{\Theta_T=\{} \mbox{ %
non-decreasing in $t$ for each $\gamma$
 }\},
\end{array}
\end{equation}
where, we define a total order $\succeq$ on the initial/boundary set 
$\Gamma_T$ by
\begin{equation}
\eqna{Gammaorder}
s\le t,\ z\le y\ \Leftrightarrow\ (0,T)\succeq (0,t)\succeq (0,s)\succeq (0,0)
\succeq (z,0)\succeq (y,0)\succeq (1,0).
\end{equation}
For example,
\begin{equation}
\eqna{Gamma_infinity}
\theta((1,0),t)=1,\ \ t\in[0,T],\ \theta\in\Theta_T.
\end{equation}

For each $\theta\in\Theta_T$, $w\in W$, and $z\in[0,1]$, define 
$\tilde{w}_{\theta,w,z}$, a non-negative valued continuous function 
of $(s,t)$ satisfying $0\le s\le t\le T$, by
\begin{equation}
\eqna{w2omega}
\tilde{w}_{\theta,w,z}(s,t)=\left\{\begin{array}{ll}\dsp w(\theta((z,0),t),t), 
& \mbox{if }\ s=0,\\ \dsp w(\theta((0,s),t),t),& \mbox{if }\ s>0.
\end{array}\right.
\end{equation}
Note that $\tilde{w}_{\theta,w,z}$ is independent of $z$ if $s>0$.
Let $\theta\in\Theta_T$, and put
\begin{equation}
\eqna{varphi_theta}
\begin{array}{l}\dsp
\varphi_{\theta}(dw,\gamma,t)=\int_{z\in[y_0,1]}
 \prb{\tilde{\nu}_{\theta,w,z}(t)
=\tilde{\nu}_{\theta,w,z}(t_0)}\,\mu_0(dw\times dz),
\\ \dsp
\ \gamma=(y_0,t_0)\in\Gamma_t,\ (\gamma,t)\in\Delta_T\,,
\end{array}
\end{equation}
where
$\dsp\tilde{\nu}_{\theta,w,z}$ is the point process
with last-arrival-time dependent intensity, denoted by $N$
in \secu{ptproc}, with the `intensity density' $\omega(s,t)$
in the definition \eqnu{intensity} of the process
given by $\tilde{w}_{\theta,w,z}(s,t)$ of \eqnu{w2omega},
and $\mu_0$ is a Borel probability measure
on the direct product space $W\times[0,1]$.

The following is proved in \cite{fluid14}.
\thmb[{\cite[Theorem~9, (89)]{fluid14}}]
\thma{FPtheorem}
There exists a unique flow $y_C\in\Theta_T$ such that 
\begin{equation}
\eqna{yCFPG}
\theta(\gamma,t)=1-\varphi_{\theta}(W,\gamma,t),
\ \gamma=(y_0,t_0)\in\Gamma,\ (\gamma,t)\in\Delta_T\,,
\end{equation}
holds for $\theta=y_C$.
\DDD
\thme
In \cite{fluid14} below the Theorem it is remarked that
there exists $\mu_t$, 
taking values in the space of Borel probability measures on $W\times[0,1]$,
such that
\begin{equation}
\eqna{varphiyC}
\arrb{l}\dsp
\varphi_{y_C}(dw,(y_0,t_0),t)=\mu_t(dw\times[y_C((y_0,t_0),t),1]),
\\ \dsp
((y_0,t_0),t)\in \Delta_T\,,
\arre
\end{equation}
and below Theorem~1 in \cite{fluid14} it is also remarked that
$\mu_t$ and $y_C$ satisfy
\begin{equation}
\eqna{yc}
y_C(\gamma,t)=y_0+\int_{t_0}^t \int_{W\times[y_C(\gamma,s),1)}
 w(z,s) \mu_s(dw\times dz)\,ds.
\end{equation}

Next we state our assumptions
on the infinite particle limit $N\to\infty$ of initial ($t=0$) conditions.
We will assume a standard supremum norm on the space of continuous
functions on the closed interval $[0,1]\times[0,T]$,
with which we define the weak convergence of Borel probability measures
on $W$ in the standard way \cite{Billingsley}.
For the initial distribution $\mu\ssN_0$ in \eqnu{initialdistributionfcn},
we assume the following weak convergence with additional
uniform bounds on the order of convergence
to $\mu_0$ in \eqnu{varphi_theta}:
\begin{equation}
\eqna{muN0conv}
\begin{array}{l}\dsp 
\mbox{For any set $H\subset C^0(W\times[0,1];\reals)$ of
uniformly bounded,}
\\ \dsp
\mbox{equicontinuous functions,}\ 
\exists \delta\in(0,\frac12),\ \exists C>0;
\\ \dsp
\ (\forall N\in\nintegers)
\ (\forall \tilde{h}\in H)
\ (\forall y\in [0,1])
\\ \dsp
\biggl|
\int_{W\times[y,1]} \tilde{h}(w,z)\,\mu\ssN_0(dw\times dy)
-\int_{W\times[y,1]} \tilde{h}(w,z)\,\mu_0(dw\times dy)\biggr|
\le \frac{C}{N^{\delta}}.
\arre
\end{equation}

Denote the marginal distribution of $\mu_0$ on $W$ by $\lambda$;
\begin{equation}
\eqna{lambda}
\lambda(dw)=\mu_0(dw\times [0,1]),
\end{equation}
and put
\begin{equation}
\eqna{lambdaN}
\lambda\ssN=\frac1N\sum_{i=1}^N\delta_{w_i}.
\end{equation}
Comparing with \eqnu{mut}, we see that $\lambda\ssN$ is the marginal
distribution of $\mu\ssN_t$ on $W$ for all $t$;
\begin{equation}
\eqna{lambdamu0N}
\lambda\ssN(dw)=\mu\ssN_{t}(dw\times [0,1]),\ t\in[0,T].
\end{equation}
For $\lambda$ we assume
\begin{equation}
\eqna{rw}
M_W:=\int_W \wnorm{w} \lambda(dw) <\infty.
\end{equation}
The assumption \eqnu{muN0conv} implies, with
\eqnu{lambdamu0N} and \eqnu{lambda},
\begin{equation}
\eqna{lambda0}
\lambda\ssN\to\lambda,\ \mbox{weakly as}.\ N\to\infty.
\end{equation}
In addition we assume convergence of the average of $\wnorm{w}$:
\begin{equation}
\eqna{lambda1}
\limf{N} \int_W\wnorm{w}\lambda\ssN(dw)=M_W\,.
\end{equation}
\remb
%
In \eqnu{muN0conv} we assume uniform order of convergence $O(N^{-\delta})$,
while Ascoli--Arzel\`a type theorem implies uniform convergence
but has no control in general on the order of convergence.
%
\DDD
\reme

We are ready to state the main results of this paper.
\thmb
\thma{HDLposdep}
Under the assumptions \eqnu{Cw}, \eqnu{muN0conv}, \eqnu{rw}, and
\eqnu{lambda1},
with probability $1$, $\mu\ssN_t\to \mu_t$, weakly as $N\to\infty$, 
uniformly in $t$, where $\mu_t$ is as in \eqnu{varphiyC}.
Explicitly, we prove 
\begin{equation}
\eqna{HDLposdepconv}
\limf{N}
\sup_{t\in[0,T]}
\biggl|
\int_W h(w) \mu\ssN_t(dw\times[y,1])-\int_W h(w) \mu_t(dw\times[y,1])
\biggr|=0,\ a.s.,
\end{equation}
for all $y\in[0,1]$ and bounded continuous function
$h:\ W\to \reals$.

Assume in addition that
for a positive integer $L$ and $y_i\in[0,1)$, $i=1,2,\ldots,L$,
\begin{equation}
\eqna{taggptclconvass}
\nu\ssN_i=\nu_i,\ N\in\nintegers,
\mbox{ and }
\limf{N} y\ssN_{i} =y_i\,,
\mbox{ for } i=1,2,\ldots,L.
\end{equation}
Then, with probability $1$, the tagged particle system 
\[
(Y\ssN_1(t),Y\ssN_2(t),\ldots,Y\ssN_L(t))
\]
converges as $N\to\infty$, uniformly in $t\in[0,T]$ to a limit
$(Y_1(t),Y_2(t),\ldots$, $Y_L(t))$. Here, for each $i=1,2,\ldots,L$, 
$Y_i$ is the unique solution to
\begin{equation}
\eqna{taggedlimit}
\arrb{l}\dsp
Y_i(t) = y_i
 + \int_{s\in(0,t]} \int_{(w,z)\in W\times [Y_i(s-).1]}
 w(z,s)\,\mu_s(dw\times dz)\,ds
\\ \dsp \phantom{Y_i(t)=}
 - \int _{s\in (0,t]} \int _{\xi \in \preals } Y_i(s-) 
\chrfcn{\xi \in [0,w_i (Y_i(s-),s))} \nu _{i}(d\xi ds),
\\ \dsp 
i=1,2,\ldots,L,\ t\in[0,T].
\arre
\end{equation}
\DDD
\thme
\thmu{HDLposdep} implies propagation of chaos for 
the stochastic ranking processes.
For each $N$ all of $\{Y_i\ssN \}$ are random and interact with each other 
and $\mu\ssN_t$ is also random.
However, the limit $\mu_t$ is deterministic.
Furthermore, the randomness of the limit process $Y_i$ of a tagged 
particle depends only on its associated Poisson random measure $\nu _i$,
and is independent of $Y_j$ or $\nu_j$ with $j\ne i$.

Incidentally, we proved almost sure convergence for $\mu\ssN_t$
in \eqnu{HDLposdepconv}, while we made no assumptions on the
relation between the set of measures $\{\nu\ssN_i\mid i=1,2,\ldots,N\}$
for different $N$. 
This is stronger than the law of large numbers appearing in the
context of random walks $\dsp \frac1N \sum_{i=1}^N X_i$, where
each $X_i$ is fixed for all $N$, while \eqnu{HDLposdepconv} corresponds to
considering $\dsp \frac1N \sum_{i=1}^N X\ssN_i$, where we assume
no relation among $X\ssN_i$ for different $N$.
Such a type of convergence is known for sums of real valued random 
variables as complete convergence \cite{HR47,Er49,Er50}.
In this context, \eqnu{varphiyC} is an example of complete convergence
for a sequence of measure valued random variables.

\section{Infinite particle limit of flow driven stochastic ranking process.}
\seca{fSRP}

\subsection{Flow driven stochastic ranking process.}

In this section, we introduce an intermediate model
which we use to prove convergence of $\mu\ssN_t$ in \thmu{HDLposdep}.

Let $\{w_i\}$, $\{y\ssN_i\}$, and $\{\nu\ssN_i\}$ be as in the
stochastic ranking process \eqnu{YNC2}.
Let $\theta\in\Theta_T$ be a flow, and 
for each $i=1,2,\ldots,N$, let $\tilde{\nu}\ssNtheta_i$ be a 
point process of $N$ in \eqnu{lastarrivaltimedepSI},
with $\nu=\nu\ssN_i$ and $\dsp\omega=\tilde{w}_{\theta,w_i,y\ssN_i}$,
where the last notation is as in \eqnu{w2omega} with 
$w=w_i$ and $z=y\ssN_i$.
Define a system of stochastic processes $Y\ssNtheta_i$, $i=1,2,\ldots,N$,
by
 \begin{equation}
\eqna{SRPYthetapre}
 \arrb{l}\dsp
 Y\ssNtheta_i(t) 
 \\ \dsp \phantom{Y}
= y\ssN_i + \frac1N\sum _{j=1} ^N
 \int _{s\in (0,t]} 
 \chrfcn{ Y\ssNtheta_j(s-)>Y\ssNtheta_i(s-) }
 \tilde{\nu}\ssNtheta_j(ds)
 \\ \dsp \phantom{Y=}
 - \int _{s\in (0,t]} Y\ssNtheta_i(s-)
 \tilde{\nu}\ssNtheta_i(ds), \ \ i=1,2,\dots ,N,\ t\geq 0.
 \arre
 \end{equation}
We will call this system, the stochastic ranking process driven by
the flow $\theta$.

The fluctuation of \eqnu{SRPYthetapre} is coupled to those of
the stochastic ranking process of \eqnu{SRPY} via the set of
Poisson random measures $\{\nu\ssN_i\}$.
Using \eqnu{arrivaltimedef} and \eqnu{lastarrivaltimedepSI},
we have an expression of $\tilde{\nu}\ssNtheta_i$ using $\nu\ssN_i$. 
Define a sequence of stopping times,
$0=\tau\ssNtheta_{i,0}<\tau\ssNtheta_{i,1}<\cdots$,
by
\begin{equation}
\eqna{latestarrivaltimeSItilde}
\arrb{l}\dsp
\tau\ssNtheta_{i,0}=0,
\\ \dsp
\tau\ssNtheta_{i,k+1}=\inf\{t>\tau\ssNtheta_{i,k}\mid 
\nu\ssN_i(
\{(s,\xi)\in(\tau\ssNtheta_{i,k},T]\times\preals\mid
\\ \dsp \phantom{\tau\ssNtheta_{i,k+1}=}
0\le \xi\le w_i(Y\ssNtheta_i(s-),s)\})>0\},\ k\in\pintegers.
\arre
\end{equation}
$\tau\ssNtheta_{i,k}$ is the time that the particle $i$ in the
flow driven stochastic ranking process jumps to the top for 
the $k$-th time. Put
 \begin{equation}
\eqna{flowstarttime}
\arrb{l}\dsp
\gamma\ssNtheta_i(t)=\left\{\arrb{ll}\dsp 
(y\ssN_i,0),& 0\le t<\tau\ssNtheta_{i,1},
\\ \dsp 
(0,\tau\ssNtheta_{i,k}),& \tau\ssNtheta_{i,k}\le t<\tau\ssNtheta_{i,k+1},
\ k=1,2,\ldots.
\arre\right.
\arre
 \end{equation}
Using $\gamma\ssNtheta_i(t)$, we have an expression
\begin{equation}
\eqna{lastarrivaltimedepSItilde}
\tilde{\nu}\ssNtheta_i(t)=\int_{s\in(0,t]} \int_{\xi\in[0,\infty)}
\chrfcn{\xi\in[0,w_i(\theta(\gamma\ssNtheta_i(s-),s-),s))}
\nu\ssN_i(d\xi\,ds),
\end{equation}
and substituting \eqnu{lastarrivaltimedepSItilde} in
\eqnu{SRPYthetapre} we further have
 \begin{equation}
\eqna{SRPYtheta}
 \arrb{l}\dsp
 Y\ssNtheta_i(t) 
= y\ssN_i
 \\ \dsp \phantom{Y=}
 + \frac1N\sum _{j=1} ^N
 \int _{s\in (0,t]} \!\int_{\xi\in[0,\infty)}
 \chrfcn{ Y\ssNtheta_j(s-)>Y\ssNtheta_i(s-) }
 \\ \dsp \phantom{
 Y= + \frac1N\sum _{j=1} ^N \int _{s\in (0,t]} \!\int_{\xi\in[0,\infty)}}\times
\chrfcn{\xi\in[0,w_j(\theta(\gamma\ssNtheta_j(s-),s-),s))}
\nu\ssN_j(d\xi\,ds)
 \\ \dsp \phantom{Y=}
 - \int _{s\in (0,t]} \!\int_{\xi\in[0,\infty)}
Y\ssNtheta_i(s-)
\chrfcn{\xi\in[0,w_i(\theta(\gamma\ssNtheta_i(s-),s-),s))}
\nu\ssN_i(d\xi\,ds),
\\ \dsp
i=1,2,\dots ,N,\ \ t\geq 0.
 \arre
 \end{equation}
This is to be compared with the (original) stochastic ranking process
\eqnu{SRPY}.
We see that \eqnu{SRPYtheta} is obtained from \eqnu{SRPY} by
replacing $Y\ssN_j(s)$ appearing as a variable for $w_j$
by $\theta(\gamma\ssNtheta_j(s),s)$, and otherwise,
by $Y\ssNtheta_j(s)$.

\subsection{Characteristic curve and distribution function.}
\seca{fSRPYC}
\seca{fSRPvarphi}

For each $i=1,2,\ldots,N$ and $0\le t_0\le t$,
define an event $J\ssNtheta_i(t_0,t)\subset\Omega$ by
\begin{equation}
\eqna{JtopnuKtheta}
J\ssNtheta_i(t_0,t)=\{ \omega\in\Omega\mid 
\tilde{\nu}\ssNtheta_i(t)(\omega)>\tilde{\nu}\ssNtheta_i(t_0)(\omega)\},
\end{equation}
Since $\tilde{\nu}\ssNtheta$ is increasing, the complement is
\begin{equation}
\eqna{JcNthetacondprobexplicit}
 J\ssNtheta_i(t_0,t)^c=\{ \omega\in\Omega\mid 
\tilde{\nu}\ssNtheta_i(t)(\omega)=\tilde{\nu}\ssNtheta_i(t_0)(\omega)\}.
\end{equation}
On the event $J\ssNtheta_i(t_0,t)^c$, the contribution of the last term
on the right hand side of \eqnu{SRPYthetapre} to the difference
$Y\ssNtheta_i(t)-Y\ssNtheta_i(t_0)$ is $0$, hence $Y\ssNtheta_i$ is
non-decreasing in the interval $[t_0,t]$.
In other words, $J\ssNtheta_i(t_0,t)$ of \eqnu{JtopnuKtheta} is the
event that the particle $i$ jumps to the top $y=0$ during $(t_0,t]$.

Define the characteristic curve
for the stochastic ranking process driven by the flow $\theta$ by
\begin{equation}
\eqna{yCNtheta}
Y\ssNtheta_C((y_0,t_0),t)
=y_0+\frac1N\sum_{j\in[1,N];\ Y\ssNtheta_j(t_0)\ge y_0}
 \chrfcn{J\ssNtheta_j(t_0,t)}\,,
\end{equation}
for $(y_0,t_0)\in [0,1]\times[0,T]$ and $t\ge t_0$.
For example,
\begin{equation}
\eqna{yCNsmp}
Y\ssNtheta_C((1,t_0),t)=1\,,\ t\ge t_0\ge 0.
\end{equation}
Since $N Y\ssNtheta_j$ takes integer values,
we can write \eqnu{yCNtheta} using \eqnu{JcNthetacondprobexplicit} as
\begin{equation}
\eqna{yCNJctheta}
Y\ssNtheta_C((y_0,t_0),t)
=y_0+\frac{[N\,(1-y_0)]}N
-\frac1N\sum_{j;\ Y\ssNtheta_j(t_0)\ge y_0}
 \chrfcn{{J\ssNtheta_j(t_0,t)}^c},
\end{equation}
where $[x]$ is the largest integer not exceeding $x$.

We note the following expression corresponding to \eqnu{YNC2}.
\lemb
\lema{YNCeqYNi}
\itmb
\item
For $t\ge t_0\ge0$ and $0\le y_0\le 1$, it holds that
\begin{equation}
\eqna{yCNxi}
\arrb{l} \dsp
Y\ssNtheta_C((y_0,t_0),t) = y_0
+ \frac1N \sum_{j=1}^N
 \int _{s\in (t_0,t]} \int _{\xi \in [0,\infty )} 
 \\ \dsp {}
\chrfcn{Y\ssNtheta_j(s-)\ge Y\ssNtheta_C((y_0,t_0),s-)} 
\chrfcn{\xi \in [0,w_j(\theta(\gamma\ssNtheta_j(s-),s),s))}
\nu\ssN _j (d\xi ds),
\arre
\end{equation}
\item
It holds that
\begin{equation}
\eqna{YNthetaieqYNCfg}
Y\ssNtheta_i(t)=Y\ssNtheta_C(\gamma\ssNtheta_i(t),t),\ t\in[0,T],
\end{equation}
where $\gamma\ssNtheta_i$ is as in \eqnu{flowstarttime}.
\DDD
\itme
\leme
\prfb
Since on the event $J\ssNtheta_i(t_0,t)^c$, 
the contribution of the last term
on the right hand side of \eqnu{SRPYthetapre} to the difference
$Y\ssNtheta_i(t)-Y\ssNtheta_i(t_0)$ disappears,
$Y\ssNtheta_i(t)-Y\ssNtheta_i(t_0)$ and
\[
Y\ssNtheta_C((Y\ssNtheta_i(t_0),t_0),t)-Y\ssNtheta_i(t_0)
\]
should
satisfy the same equation for $[t_0,t]$ on $J\ssNtheta_i(t_0,t)^c$.
This implies \eqnu{yCNxi}.

By definition \eqnu{flowstarttime}, $J\ssNtheta_i(t_0,t)^c$ holds
if $\gamma\ssNtheta_i(t)=(y_0,t_0)$. Hence \eqnu{YNthetaieqYNCfg} follows.
%
\QED
\prfe

In analogy to \eqnu{mut}, for each positive integer $N$ define
a joint empirical distribution of jump rate and position of particles
on $W\times[0,1]$ by
\begin{equation}
\eqna{distributionfcntheta}
\mu\ssNtheta_t
=\frac1N \sum_{i=1}^N \delta_{(w_i,Y\ssNtheta_i(t))}
\ \ t\in [0,T].
\end{equation}
When integrated over position, we recover the 
distribution $\lambda\ssN$ of jump rates for the (original)
stochastic ranking process \eqnu{lambdaN}, 
independently of $t$ and $\theta$;
\begin{equation}
\eqna{lambdamu0Ntheta}
\mu\ssNtheta_t(dw\times[0,1])=\lambda\ssN(dw).
\end{equation}

The initial distribution of $\mu\ssNtheta_t$,
found by substituting \eqnu{SRPYthetapre} at $t=0$ to
\eqnu{distributionfcntheta}, coincide with that of 
the original model $\mu\ssN_0$ in \eqnu{initialdistributionfcn}:
\begin{equation}
\eqna{initialdistributionfcntheta}
\mu\ssNtheta_0=\mu\ssN_0\,.
\end{equation}

A set of spatial distribution functions $\varphi\ssNtheta$ is
defined analogously to \eqnu{varphiN},
by a convolution of \eqnu{yCNtheta} and \eqnu{distributionfcntheta};
\begin{equation}
\eqna{varphiNtheta}
\arrb{l}\dsp
\varphi\ssNtheta(dw,\gamma,t)
=\mu\ssNtheta_t(dw\times[Y\ssNtheta_C(\gamma,t),1]),
\\ \dsp
\gamma=(y_0,t_0)\in [0,1]\times[0,T],\ t\ge t_0.
\arre
\end{equation}
$\varphi\ssNtheta(dw,\gamma,t)$
denotes the empirical distribution of jump rates (intensity densities)
of those particles which was in a downstream of $y_0$ at time $t_0$
that have not jumped to the top in the time period $(t_0,t]$.
For $t=t_0$ we have
\begin{equation}
\eqna{varphiNt0theta}
\varphi\ssNtheta(dw,(y_0,t_0),t_0)=\mu\ssNtheta_{t_0}(dw\times[y_0,1]),
\ \ (y_0,t_0)\in [0,1]\times[0,T].
\end{equation}

$\varphi\ssNtheta$ is a refinement of $Y\ssNtheta_C$ in the 
following sense.
First, \eqnu{varphiNtheta} and \eqnu{distributionfcntheta} imply
\[ 
\varphi\ssNtheta(dw,\gamma,t)=
\frac1N\sum_{j=1}^N \chrfcn{Y\ssNtheta_j(t)\ge Y\ssNtheta_C(\gamma,t)}
\,\delta_{w_j}(dw).
\]
Next,
$Y\ssNtheta_j(t)\ge Y\ssNtheta_C(\gamma,t)$ occurs if and only if
the inequality holds at $t=t_0$ and  $j$ does not jump to the top in
the interval $(t_0,t]$. Hence
\eqnu{JtopnuKtheta} implies
\begin{equation}
\eqna{varphiNthetaJ}
\varphi\ssNtheta(dw,(y_0,t_0),t)
=\frac1N\sum_{j;\ Y\ssNtheta_j(t_0)\ge y_0}
 \chrfcn{J\ssNtheta_j(t_0,t)^c}\,\delta_{w_j}(dw).
\end{equation}
Comparing 
\eqnu{varphiNthetaJ} with \eqnu{yCNJctheta}, we have
\begin{equation}
\eqna{varphiNYNCtheta}
\arrb{l}\dsp
Y\ssNtheta_C(\gamma,t)=y_0+\frac{[N(1-y_0)]}{N}
 -\varphi\ssNtheta(W,\gamma,t),
\\ \dsp
\gamma=(y_0,t_0)\in[0,1]\times[0,T],\ t\ge t_0.
\arre
\end{equation}

The spatial distribution function $\varphi\ssNtheta(dw,\gamma,t)$ is defined
for any $\gamma=(y_0,t_0)\in [0,1]\times[0,T]$ satisfying $t\ge t_0$,
but is particularly important when $\gamma\in\Gamma_t$.
In fact, both for the case $\gamma=(y_0,0)\in\Gamma_i$ and the case
$\gamma=(0,t_0)\in\Gamma_t$ we have
\begin{equation}
\eqna{varphiNthetaJGamma}
\arrb{l}\dsp
\varphi\ssNtheta(dw,(y_0,t_0),t)
=\frac1N\sum_{j;\ y\ssN_j\ge y_0}
 \chrfcn{J\ssNtheta_j(t_0,t)^c}\,\delta_{w_j}(dw),
\ \ 
\\ \dsp
\gamma=(y_0,t_0)\in\Gamma_t,\ t\in[0,T].
\arre
\end{equation}
Note the conditions of the summation, which is 
non-random in \eqnu{varphiNthetaJGamma},
while is random in \eqnu{varphiNthetaJ}.
Since $J\ssNtheta_j(t_0,t)$ is independent of $\nu\ssN_i$, $i\ne j$,
\eqnu{varphiNthetaJGamma} implies that
\textit{$\varphi\ssNtheta(dw,\gamma,t)$ is an arithmetic average of
independent random variables, if $\gamma\in\Gamma_t$}.
With this fact,
we restrict the domain of definition to $\Delta_T$ of \eqnu{domyC}
in \secu{fSRPvarphiLLN}.

\subsection{Convergence of distribution function.}
\seca{fSRPvarphiLLN}

\subsubsection{Statement of the theorem.}

Here we will state a strong law of large numbers 
for the spatial distribution function $\varphi\ssNtheta$ of the
flow driven stochastic ranking process,
uniform in initial point $\gamma$ and time $t$.

For a bounded continuous function $h:\ W\to\reals$, we put
\begin{equation}
\eqna{hdata}
C_h=\sup_{w\in W} |h(w)|<\infty,
\end{equation}
and use a notation
\begin{equation}
\eqna{varphiNthetah}
\varphi\ssNtheta(h,\gamma,t)=\int_W h(w)\,\varphi\ssNtheta(dw,\gamma,t).
\end{equation}
With \eqnu{varphiNthetaJ} we have
\begin{equation}
\eqna{varphiNthetahexplicit}
\varphi\ssNtheta(h,\gamma,t)
=\frac1N\sum_{j;\ Y\ssNtheta_j(t_0)\ge y_0} h(w_j)\,
 \chrfcn{J\ssNtheta_j(t_0,t)^c},\ \ \gamma=(y_0,t_0).
\end{equation}
As in \eqnu{varphiNthetah} we also use a notation
\begin{equation}
\eqna{varphithetah}
\varphi_{\theta}(h,\gamma,t)=\int_W h(w)\,\varphi_{\theta}(dw,\gamma,t).
\end{equation}
for $\varphi_{\theta}$ in \eqnu{varphi_theta}.
\thmb
\thma{fSRPLLN}
Assume \eqnu{Cw}, \eqnu{muN0conv}, \eqnu{rw}, and \eqnu{lambda1}. 
Then for any $p>0$ there exists a positive constant $C$
depending only on $p$ and $\delta$, 
(and is independent of $N$, $\theta$, and $h$,)
such that for any bounded continuous $h:\ W\to\reals$,
\begin{equation}
\eqna{fSRPHDLepsilon}
\EE{ \sup_{(\gamma,t)\in\Delta_T}\biggl|
\varphi\ssNtheta(h,\gamma,t)-\varphi_{\theta}(h,\gamma,t)
\biggr|^{2p} }
\le \frac{C\,C_h^{2p}}{N^{2p\delta}}\,,
\ N\in\nintegers,
\end{equation}
holds, where $C_h$ is as in \eqnu{hdata}.
\DDD
\thme
\thmu{FPtheorem}, \thmu{fSRPLLN} with $h=\chrfcn{W}$ and $\theta=y_C$, and 
\eqnu{varphiNYNCtheta} imply the following.
\corb
\cora{fSRPyCHDLepsilon}
Under the assumptions of \thmu{fSRPLLN},
for any $p>0$
there exists a positive constant $C$ depending only on $p$ and $\delta$,
such that
\begin{equation}
\eqna{fSRPyCHDLepsilon}
\EE{ \sup_{(\gamma,t)\in\Delta_T}\biggl|
Y_C\ssNyC(\gamma,t)-y_C(\gamma,t)
\biggr|^{2p} }
\le \frac{C}{N^{2p\delta}}\,,
\ N\in\nintegers.
\end{equation}
Among $\theta\in\Theta$, $\theta=y_C$ is the only flow that
satisfies \eqnu{fSRPyCHDLepsilon}.
\DDD
\core
As in the proof of the main theorem \thmu{HDLposdep},
\thmu{fSRPLLN} implies convergence of joint empirical distribution
$\mu\ssNtheta_t$ to $\mu_{y_C,t}$, a deterministic distribution 
determined by $\varphi_{y_C}(h,\gamma,t)$.
For $\theta\ne y_C$, the stochastic ranking process driven by the
flow $\theta$ converges to a limit with deterministic distribution, but
the resulting trajectories of particles are different from the given flow
$\theta$, due to the uniqueness result of \thmu{FPtheorem}.

A proof of \thmu{fSRPLLN} is composed of $2$ parts.
In \secu{fSRPLLN} we prove that 
$\varphi\ssNtheta-\EE{\varphi\ssNtheta}$ converges to $0$,
using a strong uniform law of large numbers in \cite{LLN15}.
In \secu{fSRPvarphiexpec} we prove that 
$\EE{\varphi\ssNtheta}$ converges to $\varphi_{\theta}$,
using the estimates in \cite{fluid14}.
Relevant results of \cite{fluid14} are summarized in \secu{ptproc}
for convenience.

\subsubsection{Strong uniform law of large numbers.}
\seca{fSRPLLN}

Here we will prove the following.
\thmb
\thma{indepPoissonLLNepsilon}
Assume \eqnu{Cw}, \eqnu{muN0conv}, \eqnu{rw}, and \eqnu{lambda1}. Then
for any $p$ and $\delta$ satisfying
\begin{equation}
\eqna{indepPoissonLLNepsilon_cond}
p>0,\ \ 0<\delta<\frac12\,,
\end{equation}
there exists a positive constant $C$
(independent of $N$, $\theta$, and $h$,)
such that for any bounded continuous $h:\ W\to\reals$
\begin{equation}
\eqna{indepPoissonLLNepsilon}
\EE{ \sup_{(\gamma,t)\in\Delta_T}\biggl|
\varphi\ssNtheta(h,\gamma,t)-\EE{\varphi\ssNtheta(h,\gamma,t)}
\biggr|^{2p} }
\le \frac{C\,C_h^{2p}}{N^{2p\delta}}\,,
\ N\in\nintegers,
\end{equation}
where $C_h$ is as in \eqnu{hdata}.
\DDD
\thme
\prfb
Put $\dsp h_{\pm}=(\pm h)\vee 0$, so that $h=h_+-h_-$ decomposes the function
$h$ to positive and negative parts.
Applying \eqnu{varphiNthetaJ}, \eqnu{domyC}, \eqnu{Gammat},
and triangular inequality in the form
\[
(a+b)^{q}\le (2^{q-1}\vee 1)(a^{q}+b^{q}),
\ a\ge0,\ b\ge0,\ q>0,
\]
we have, with $\gamma=(y_0,t_0)\in\Gamma_t$,
\begin{equation}
\eqna{FISRPLLN1}
\arrb{l}\dsp
\EE{ \sup_{((y_0,t_0),t)\in\Delta_T}\biggl|
\varphi\ssNtheta(h,\gamma,t)-\EE{\varphi\ssNtheta(h,\gamma,t)}
\biggr|^{q} }
\\ \dsp {}
=\EE{ \sup_{((y_0,t_0),t)}\biggl|\frac1N
 \sum_{Y\ssNtheta_j(t_0)\ge y_0} (h_+(w_j)-h_-(w_j))\,
\\ \dsp \phantom{=\EE{ \sup_{((y_0,t_0),t)}\biggl|\frac1N \sum}}\times
( \chrfcn{J\ssNtheta_j(t_0,t)^c}-\EE{\chrfcn{J\ssNtheta_j(t_0,t)^c}})
\biggr|^{q} }
\\ \dsp {}
\le (2^{q-1}\vee1)\EE{ \sup_{((y_0,t_0),t)}\biggl|\frac1N
 \sum_{Y\ssNtheta_j(t_0)\ge y_0} h_+(w_j)\,
\\ \dsp \phantom{
  \le (2^{q-1}\vee1)\EE{ \sup_{((y_0,t_0),t)}\biggl|\frac1N \sum}}\times
( \chrfcn{J\ssNtheta_j(t_0,t)^c}-\EE{\chrfcn{J\ssNtheta_j(t_0,t)^c}})
\biggr|^{q} }
\\ \dsp \phantom{\le}
+ (2^{q-1}\vee1)\EE{ \sup_{((y_0,t_0),t)}\biggl|\frac1N
 \sum_{Y\ssNtheta_j(t_0)\ge y_0} h_-(w_j)\,
\\ \dsp \phantom{
  \le (2^{q-1}\vee1)\EE{ \sup_{((y_0,t_0),t)}\biggl|\frac1N \sum}}\times
( \chrfcn{J\ssNtheta_j(t_0,t)^c}-\EE{\chrfcn{J\ssNtheta_j(t_0,t)^c}})
\biggr|^{q} }
\\ \dsp {}
\le (2^{q-1}\vee1)
\,(R_{q,1,+}\ssN+R_{q,1,-}\ssN+R_{q,2,+}\ssN+R_{q,2,-}\ssN),
\arre
\end{equation}
where
\begin{equation}
\eqna{FISRPLLNR1}
R_{q,1,\pm}\ssN=\EE{\sup_{t\in[0,T]}\sup_{0\le y_0<1}\biggl|\frac1N
 \sum_{i;\ y\ssN_i\ge y_0} h_{\pm}(w_i)\,
( \chrfcn{J\ssNtheta_i(0,t)^c}-\EE{\chrfcn{J\ssNtheta_i(0,t)^c}})
\biggr|^{q} }
\end{equation}
and
\begin{equation}
\eqna{FISRPLLNR2}
R_{q,2,\pm}\ssN=\EE{\sup_{t\in[0,T]}\sup_{0\le t_0<t}\biggl|\frac1N
 \sum_{i=1}^N h_{\pm}(w_i)\,
( \chrfcn{J\ssNtheta_i(t_0,t)^c}-\EE{\chrfcn{J\ssNtheta_i(t_0,t)^c}})
\biggr|^{q} }.
\end{equation}

To bound \eqnu{FISRPLLNR1} and \eqnu{FISRPLLNR2}, we refer to
the last theorem in \cite[\S2]{LLN15}.
We reproduce the theorem in a specific form of
\[
Z\ssN_i(s,t)=a\ssN_i\chrfcn{\nu\ssN_i(t)>\nu\ssN_i(s)}
\]
in place of $Z\ssN_i(s,t)$ in the reference.
\prpb[{\protect\cite[Theorem~3]{LLN15}}]
\prpa{LLN15r}
Let $T>0$, and 
for each $N\in\nintegers$, let $\nu\ssN_i$, $i=1.\ldots,N$, be 
a sequence of independent random variables taking values in a space of
non-negative valued non-decreasing right continuous functions on $[0,T]$
with left limit.
Let $r>0$, and for $N\in\nintegers$, let $M\ssN>0$ and let
$a\ssN_i$ and $w\ssN_i$, $i=1,\ldots,N$, be non-negative sequences.
Assume that $a\ssN_i\le M\ssN$, $i=1,\ldots,N$, and that
\[\arrb{l}\dsp
|\prb{\nu\ssN_i(t_2)>\nu\ssN_i(t_1)}-
\prb{\nu\ssN_i(s_2)>\nu\ssN_i(s_1)}|
\\ \dsp
\le w\ssN_i(|t_1-s_1|^r+|t_2-s_2|^r),
\arre\]
for $0\le t_1\le t_2\le T$ and $0\le s_1\le s_2\le T$.
Then for any $\dsp \delta\in(0,\frac12)$ and $p>0$,
\begin{equation}
\eqna{main2}
\arrb{l}\dsp
\EE{ \sup_{t_2\in[0,T]}\sup_{0\le t_1<t_2}
\biggl|\frac1N\sum_{i=1}^N
a\ssN_i\, 
(\chrfcn{\nu\ssN_i(t_2)>\nu\ssN_i(t_1)}
-\prb{\nu\ssN_i(t_2)>\nu\ssN_i(t_1)})\biggr|^{p}}^{1/p}
\\ \dsp {}
\le \frac{M\ssN}{N^{\delta}}2^{1-1/q}
 (C_q^q\,(2T(\bar{w}\ssN)^{1/r}+1)+2^{2q})^{1/q},
\\ \dsp {}
N=N_0,N_0+1,\ldots,
\arre
\end{equation}
where 
$\dsp q=q(p,\delta)= 3\vee \frac{r+1}r \,\frac{2\delta}{1-2\delta} \vee p$,
$N_0$ is the smallest integer
satisfying $N_0^{rq/(2rq+2r+2)}\ge 2$, 
$\dsp\bar{w}\ssN=\frac1N\sum_{i=1}^N w\ssN_i$,
and
$\dsp C_q=\biggl(\frac12(4k)^q+\frac{2k}{2k-q}(8k)^q\biggr)^{1/q}$
with $k$ the smallest integer greater than $\dsp \frac12q$,
(In particular, $q$ and $N_0$ are independent of 
$N$, $M\ssN$, and $\bar{w}\ssN$.)
\DDD
\prpe

To bound $R_{q,2,\pm}\ssN$ in \eqnu{FISRPLLNR2},
we apply \prpu{LLN15r} with $\nu\ssN_i=\tilde{\nu}\ssNtheta_i$
and $\dsp a\ssN_i=h_{\pm}(w_i)$.
Note that
\eqnu{st_dep_diff} and \eqnu{w2omega} with $w=w_i$ and $z=y\ssN_i$, 
and \eqnu{Wnorm} imply
\begin{equation}
\eqna{indepPoissonLLNepsilonprf2}
\arrb{l}\dsp
0\le -\pderiv{}{t}\prb{\tilde{\nu}\ssNtheta_i(t)=\tilde{\nu}\ssNtheta_i(s)} 
\le  \norm{\tilde{w}_{\theta,w_i,y\ssN_i}}\le \wnorm{w_i},
\\ \dsp
0\le \pderiv{}{s}\prb{\tilde{\nu}\ssNtheta_i(t)=\tilde{\nu}\ssNtheta_i(s)} 
\le  \norm{\tilde{w}_{\theta,w_i,y\ssN_i}}\le \wnorm{w_i}.
\arre
\end{equation}
Recall \eqnu{JtopnuKtheta} and \eqnu{JcNthetacondprobexplicit}.
Comparing the left hand side of \eqnu{main2} with the right hand side of
\eqnu{FISRPLLNR2} with $q=2p$,
we see that we can apply \prpu{LLN15r} to $R_{2p,2,\pm}\ssN$
with
\begin{equation}
\eqna{indepPoissonLLNepsilonprf6}
M\ssN=C_h,\ \ r=1,\ \ w\ssN_i=\wnorm{w_i}\,.
\end{equation}
\prpu{LLN15r} then implies that
for any $\dsp \delta\in(0,\frac12)$ and $p>0$,
\begin{equation}
\eqna{indepPoissonLLNepsilonprf4}
\arrb{l}\dsp
R_{2p,2,\pm}\ssN\le 
\frac{C_h^{2p}}{N^{2p\delta}}2^{2p-1} \biggl(
C_{q}^{q}(2T\bar{w}{}\ssN+1)+2^{2q}
\biggr)^{2p/q},
\ \ N> 2^{2+(4/q)},
\arre
\end{equation}
where 
\[
\arrb{l}\dsp
q=q(p,\delta)=3\vee \frac{4\delta}{1-2\delta} \vee (2p),
\\ \dsp
\bar{w}{}\ssN
=\frac1N \sum_{i=1}^N \wnorm{w_i}=\int_W\wnorm{w}\lambda\ssN(dw),
\arre
\]
and $C_{q}$ is a positive constant depending only on $q$.

Combining
\eqnu{indepPoissonLLNepsilonprf6}C
\eqnu{indepPoissonLLNepsilonprf4}, and
\eqnu{lambda1},
we see that there exists a positive constant $C_{p,\delta}$
independent of $N$, $\theta$, and $h$, such that
\begin{equation}
\eqna{indepPoissonLLNepsilonprf9}
R_{2p,2,\pm}\ssN\le \frac{C_{p,\delta}\,C_h^{2p}}{N^{2p\delta}},
\ N\in\nintegers.
\end{equation}

To bound $R_{q,1,\pm}\ssN$ in \eqnu{FISRPLLNR1}, we first note that
\begin{equation}
\eqna{indepPoissonLLNepsilonprfn1}
\arrb{l}\dsp
R_{q,1,\pm}\ssN
\\ \dsp
\le\sum_{j=1}^N\EE{\sup_{t\in[0,T]}\biggl|\frac1N
 \sum_{i;\ N\,y\ssN_i\ge j-1} h_{\pm}(w_i)\,
( \chrfcn{J\ssNtheta_i(0,t)^c}-\EE{\chrfcn{J\ssNtheta_i(0,t)^c}})
\biggr|^{q} }
\\ \dsp
\le\sum_{j=1}^N\EE{\sup_{t\in[0,T]}\sup_{0\le t_0<t}\biggl|\frac1N
 \sum_{i;\ N\,y\ssN_i\ge j-1} h_{\pm}(w_i)\,
( \chrfcn{J\ssNtheta_i(t_0,t)^c}-\EE{\chrfcn{J\ssNtheta_i(t_0,t)^c}})
\biggr|^{q} }.
\arre
\end{equation}
Comparing \eqnu{indepPoissonLLNepsilonprfn1} with \eqnu{FISRPLLNR2},
we see that we can apply \prpu{LLN15r}
with $\nu\ssN_i=\tilde{\nu}\ssNtheta_i$
and 
\[
a\ssN_i=h_{\pm}(w_i)\chrfcn{N\,y\ssN_i\ge j-1}\,,
\]
and \eqnu{indepPoissonLLNepsilonprf6},
to the $j$-th term of the summation in the right hand side of
\eqnu{indepPoissonLLNepsilonprfn1},
in a similar way as we did to $R_{2p,2,\pm}\ssN$.
Using monotonicity of $L^p$ norms with respect to $p$
before applying \prpu{LLN15r}, we have, for $q_0\ge 2p$ and 
$\dsp \delta'\in(0,\frac12)$,
\[
R_{2p,1,\pm}\le R_{q_0,1,\pm}^{2p/q_0}
\le \frac{(C_{q_0/2,\delta'}\,C_h^{q_0})^{2p/q_0}}{N^{2p\delta'-2p/q_0}},
\ N\in\nintegers,
\]
in place of \eqnu{indepPoissonLLNepsilonprf9}.
(The extra factor $N$ compared to \eqnu{indepPoissonLLNepsilonprf9}
is from the summation with respect to $j$ 
in \eqnu{indepPoissonLLNepsilonprfn1}.)
Now choose $\delta'$ and $q_0$ to satisfy $\dsp \delta<\delta'<\frac12$
and $\dsp 0<\frac1{q_0}<\delta'-\delta$ to find
\begin{equation}
\eqna{indepPoissonLLNepsilonprf1}
R_{2p,1,\pm}\le R_{q_0,1,\pm}^{2p/q_0}
\le \frac{(C_{q_0/2,\delta'}\,C_h^{q_0})^{2p/q_0}}{N^{2p\delta}},
\ N\in\nintegers.
\end{equation}

\thmu{indepPoissonLLNepsilon} finally follows from
\eqnu{FISRPLLN1}, \eqnu{indepPoissonLLNepsilonprf1}, and
\eqnu{indepPoissonLLNepsilonprf9}.
\QED
\prfe

\subsubsection{Convergence of expectation.}
\seca{fSRPvarphiexpec}

Here we complete a proof of \thmu{fSRPLLN}.

\lemb
\lema{ufmintgbl}
Assume \eqnu{lambda1}. Then,
if $h: W\to\reals$ is bounded and continuous,
\begin{equation}
\eqna{lambda3}
\limf{N} \int_Wh(w)\,\wnorm{w}\lambda\ssN(dw)=
\int_Wh(w)\,\wnorm{w}\lambda(dw).
\end{equation}
\DDD
\leme
\prfb
Note first that for $M>0$ 
\begin{equation}
\eqna{veeplus}
a=(a\wedge M) +(a-M)_+
\end{equation}
holds.
Since
$\dsp w\mapsto \wnorm{w}\wedge M$ is bounded and continuous,
\eqnu{lambda1} (weak convergence) implies
\begin{equation}
\eqna{ufmintgblprf1}
\limf{N}\int_{W} \wnorm{w}\wedge M\,\lambda\ssN(dw)
=\int_{W} \wnorm{w}\wedge M\,\lambda(dw),
\end{equation}
which, with convergence of expectation in \eqnu{lambda1}, 
further implies
\begin{equation}
\eqna{ufmintgblprf2}
\limf{N}\int_{W} (\wnorm{w}-M)_+\,\lambda\ssN(dw)
=\int_{W} (\wnorm{w}-M)_+\,\lambda(dw).
\end{equation}
On the other hand,
dominated convergence theorem,
$\dsp 0\le (\wnorm{w}-M)_+\le \wnorm{w}$, and
\eqnu{rw} imply
\[
\limf{M} \int_W  (\wnorm{w}-M)_+ \lambda(dw)=0.
\]
Hence, for any $\eps>0$ there exists $M>0$ such that
$\dsp
\int_W  (\wnorm{w}-M)_+ \lambda(dw)<\eps. 
$
This and \eqnu{ufmintgblprf2} further imply that there exists
$N_0>0$ such that
\[
\int_W  (\wnorm{w}-M)_+ \lambda\ssN(dw)<2\eps,\ \ N\ge N_0\,.
\]
Put $\dsp C=\sup_{w\in W} |h(w)|<\infty$. Then
\[\arrb{l}\dsp
|\int_Wh(w)\wnorm{w}\lambda\ssN(dw)-\int_Wh(w)\wnorm{w}\lambda(dw)|
\\ \dsp \phantom{}
\le |\int_{W} \wnorm{w}\wedge M\,\lambda\ssN(dw)
-\int_{W} \wnorm{w}\wedge M\,\lambda(dw)|
\\ \dsp \phantom{\le}
+C\int_W  (\wnorm{w}-M)_+ \lambda\ssN(dw)
+C\int_W  (\wnorm{w}-M)_+ \lambda(dw)
\\ \dsp \phantom{}
\le |\int_{W} \wnorm{w}\wedge M\,\lambda\ssN(dw)
-\int_{W} \wnorm{w}\wedge M\,\lambda(dw)|
+3C\eps.
\arre\]
This and \eqnu{ufmintgblprf1} imply
\[
\limsup_{N\to\infty}
|\int_Wh(w)\wnorm{w}\lambda\ssN(dw)-\int_Wh(w)\wnorm{w}\lambda(dw)|
\le 3C\eps.
\]
Since the left hand side is independent of $N$, $M$, and $\eps$, 
it must be $0$.
\QED
\prfe
\prfofb{\protect\thmu{fSRPLLN}}
Let $h:\ W\to\reals$ be a bounded continuous function,
and put $\dsp C_h=\sup_{w\in W} |h(w)|$ as in \eqnu{hdata}.
\thmu{indepPoissonLLNepsilon} implies that 
to prove \thmu{fSRPLLN}, it suffices to prove
\begin{equation}
\eqna{fSRPvarphiexpec_conv}
\sup_{(\gamma,t)\in\Delta_T}\biggl|
\EE{\varphi\ssNtheta(h,\gamma,t)}-\varphi_{\theta}(h,\gamma,t)
\biggr|
\le \frac{C\,C_h}{N^{\delta}}\,,
\ N\in\nintegers,
\end{equation}
for $C$ depending only on $p$ and $\delta$
(independent of $N$, $\theta$, and $h$).

Comparing 
the definition of $\tilde{\nu}\ssNtheta_i$ given above \eqnu{SRPYthetapre}
with that of $\tilde{\nu}_{\theta,w,z}$ given below \eqnu{varphi_theta},
we see that $\tilde{\nu}\ssNtheta_i$ and $\tilde{\nu}_{\theta,w_i,y\ssN_i}$
have identical distribution.
Therefore, using 
\eqnu{JcNthetacondprobexplicit} and \eqnu{initialdistributionfcn}
in \eqnu{varphiNthetahexplicit}, we have, for $t\ge t_0$,
\begin{equation}
\eqna{varphiNthetahexpecdef}
\arrb{l}\dsp
\EE{\varphi\ssNtheta(h,(y_0,t_0),t)}
=\int_{W\times [y_0,1]} 
h(w)\, \prb{\tilde{\nu}_{\theta,w,z}(t)=\tilde{\nu}_{\theta,w,z}(t_0) }
\mu\ssN_0(dw \times dz).
\arre
\end{equation}

Let $\dsp\gamma=(y_0,t_0)\in\Gamma$ and $\dsp (\gamma,t)\in\Delta_T$.
\eqnu{varphiNthetahexpecdef} and \eqnu{varphi_theta} imply
\[\arrb{l}\dsp
\EE{\varphi\ssNtheta(h,\gamma,t)}-\varphi_{\theta}(h,\gamma,t)
\\ \dsp {}
=\int_{W\times [y_0,1]} 
h(w)\, \prb{\tilde{\nu}_{\theta,w,z}(t)=\tilde{\nu}_{\theta,w,z}(t_0) }
\,(\mu\ssN_0(dw \times dz)-\mu_0(dw\times dz)).
\arre\]
Hence,
\begin{equation}
\eqna{varphithetaexpeclimitdiffbd}
\arrb{l}\dsp
\sup_{(\gamma,t)\in\Delta_T}|\EE{\varphi\ssNtheta(h,\gamma,t)}-
\varphi_{\theta}(h,\gamma,t)|
\\ \dsp {}
\le
\sup_{((y_0,t_0),t)\in\Delta_T}
\biggl|\int_{W\times[y_0,1]}\tilde{h}_{t_0,t}(w,z)
(\mu\ssN_0(dw\times dz)-\mu_0(dw\times dz)\biggr|,
\arre
\end{equation}
where
\begin{equation}
\eqna{contilimsoln}
\tilde{h}_{t_0,t}(w,z)=
h(w)\, \prb{\tilde{\nu}_{\theta,w,z}(t)=\tilde{\nu}_{\theta,w,z}(t_0) }.
\end{equation}
Choose the set $H$ in \eqnu{muN0conv} as the set of the functions
$\tilde{h}_{t_0,t}$ in \eqnu{contilimsoln}:
\begin{equation}
\eqna{H}
H=\{\tilde{h}_{t_0,t}:\ W\times[0,1]\to\reals \mid 0\le t_0\le t\le T\}.
\end{equation}
Uniform boundedness of the functions in $H$ is obvious. 
If we prove that $H$ is also equicontinuous,
then by the assumption of \thmu{fSRPLLN} 
the consequence of \eqnu{muN0conv} holds,
which implies
\[\arrb{l}\dsp
\sup_{y_0\in[0,1]}\sup_{0\le t_0\le t\le T}
\biggl|
\int_{W\times[y_0,1]}\tilde{h}_{t_0,t}(w,z)
(\mu\ssN_0(dw\times dz)-\mu_0(dw\times dz)\biggr|
\le \frac{C\,C_h}{N^{\delta}},\ \ 
\\ \dsp
N\in\nintegers.
\arre\]
Applying this estimate to \eqnu{varphithetaexpeclimitdiffbd},
we have \eqnu{fSRPvarphiexpec_conv},
which proves \thmu{fSRPLLN}.

We are left with proving equicontinuity of $H$. 

First,
for $(w,z)\in W\times[0,1]$ and
$\dsp\tilde{w}_{\theta,w,z}$ as in \eqnu{w2omega}
(i.e., the `intensity density' for $\dsp\tilde{\nu}_{\theta,w,z}$),
and $0\le s\le t\le T$, put
\begin{equation}
\eqna{omegaOmega}
\Omega_{\theta,w,z}(s,t)=\int_s^t \tilde{w}_{\theta,w,z}(s,u)\,du
\ \mbox{ and }\ 
\tilde{\Omega}_w(s,t)=\int_s^t w(1,u)\,du.
\end{equation}
Then \eqnu{Cw} and a mean value theorem imply
\begin{equation}
\eqna{equicontiprf1}
w(1,t)-C_W\le \tilde{w}_{\theta,w,z}(s,t) \le w(1,t)+C_W,
\end{equation}
and
\begin{equation}
\eqna{equicontiprf2}
e^{-\tilde{\Omega}_w(s,t)-C_W\,(t-s)}
\le e^{-\Omega_{\theta,w,z}(s,t)}
\le e^{-\tilde{\Omega}_w(s,t)+C_W\,(t-s)}.
\end{equation}
Note also an elementary formula in \cite[(53)]{fluid14}
\begin{equation}
\eqna{equicontiprf3}
\int_{0\le u_1\le u_2\le \cdots\le u_k\le s}
\prod_{i=1}^k f(u_i)du_1\,du_2\ldots du_k
=\frac1{k!}\biggl(\int_0^s f(v)dv\biggr)^k,
\end{equation}
valid for any integrable function $f:\ \reals\to\reals$,
$s\ge0$, and $k=1,2,\ldots$,

A proof of equicontinuity of $H$ now goes in a similar way as 
that of \cite[Lemma~12]{fluid14}.
Applying \eqnu{numberdepPoissondecayposdep} to \eqnu{contilimsoln},
we have
\begin{equation}
\eqna{equicontiprf4}
\arrb{l}\dsp
\tilde{h}_{t_0,t}(w,z)
=h(w)
\sum_{k\ge0} \int_{0=:u_k< u_{k-1}<u_{k-2}<\cdots<u_1<u_0\le t_0}
\\ \dsp \phantom{\tilde{h}_{t_0,t}(w,z)=}\times
e^{-\sum_{i=0}^{k-1} \Omega_{\theta,w,z}(u_{i+1},u_i)
 -\Omega_{\theta,w,z}(u_0,t)}\,
\biggl(\prod_{i=0}^{k-1} \tilde{w}_{\theta,w,z}(u_{i+1},u_i)\,du_i\biggr).
\arre
\end{equation}
Using \eqnu{equicontiprf1}, \eqnu{equicontiprf2},
and \eqnu{omegaOmega} to \eqnu{equicontiprf4},
while noting that \eqnu{w2omega} implies that
$\tilde{w}_{\theta,w,z}(s,t)$ and $\Omega_{\theta,w,z}(s,t)$ 
is independent of $z$ if $s>0$,
we have
\[
|\tilde{h}_{t_0,t}(w,z')-\tilde{h}_{t_0,t}(w,z)|\le 
I_{11}(z,z')+I_{12}(z,z'),
\]
where
\[
\arrb{l}\dsp
I_{11}(z,z')=C_h
e^{-\tilde{\Omega}_w(0,t)+C_Wt}
\sum_{k\ge1} \int_{0=:u_k< u_{k-1}<u_{k-2}<\cdots<u_1<u_0\le t_0}
\\ \dsp \phantom{I_{11}(z,z')=} \times
\biggl(\prod_{i=0}^{k-2} (w(1,u_i)+C_W)du_i\biggr)
\\ \dsp \phantom{I_{11}(z,z')=} \times
|\tilde{w}_{\theta,w,z'}(0,u_{k-1})-\tilde{w}_{\theta,w,z}(0,u_{k-1})|
du_{k-1}\,,
\arre
\]
and
\[
\arrb{l}\dsp
I_{12}(z,z')=C_h
\sum_{k\ge0} \int_{0=:u_k< u_{k-1}<u_{k-2}<\cdots<u_1<u_0\le t_0}
e^{-\tilde{\Omega}_w(u_{k-1},t)+C_W(t-u_{k-1})}
\\ \dsp \phantom{I_{11}(z,z')=} \times
\biggl(\prod_{i=0}^{k-1} (w(1,u_i)+C_W)du_i\biggr)
\\ \dsp \phantom{I_{11}(z,z')=} \times
|e^{-\Omega_{\theta,w,z'}(0,u_{k-1})}-e^{-\Omega_{\theta,w,z}(0,u_{k-1})}|
du_{k-1}\,.
\arre
\]
Using \eqnu{w2omega}, \eqnu{Cw},
\eqnu{equicontiprf3}, 
and \eqnu{omegaOmega}, we have
\[
\arrb{l}\dsp
I_{11}(z,z')\le 
C_hC_W e^{-\tilde{\Omega}_w(t_0,t)+C_W(t+t_0)}\int_0^{t_0}
|\theta((z',0),v)-\theta((z,0),v)|dv
\\ \dsp \phantom{I_{11}(z,z')}
\le C_hC_W e^{2C_WT}\int_0^T|\theta((z',0),v)-\theta((z,0),v)|dv.
\arre
\]
Using in addition 
\begin{equation}
\eqna{equicontiprf5}
\begin{array}{l} \dsp
|e^{-x'}-e^{-x}|=e^{-(x'\wedge x)}-e^{-(x'\vee x)}
= e^{-(x'\wedge x)} \, (1- e^{-|x'-x|})
\\ \dsp \le e^{-(x'\wedge x)} |x'-x| \le e^{-x} e^{|x'-x|} |x'-x|,
\end{array}
\end{equation}
which follows from
$\dsp |x'-x|=(x'\vee x)-(x'\wedge x)\ge x-(x\wedge x')$,
we similarly have
\[
\arrb{l}\dsp
I_{12}(z,z')\le 
C_hC_W e^{-\tilde{\Omega}_w(t_0,t)+C_W(t+t_0)}
\int_0^{t_0}|\theta((z',0),v)-\theta((z,0),v)|dv
\\ \dsp \phantom{I_{12}(z,z')\le}\times 
e^{C_W\int_0^{t_0}|\theta((z',0),v)-\theta((z,0),v)|dv}
\\ \dsp \phantom{I_{12}(z,z')}
\le C_hC_W e^{2C_WT}\int_0^T|\theta((z',0),v)-\theta((z,0),v)|dv
\\ \dsp \phantom{I_{12}(z,z')\le}\times 
e^{C_W\int_0^{T}|\theta((z',0),v)-\theta((z,0),v)|dv}.
\arre
\]
Since the right hand sides of the bounds for $I_{11}$ and $I_{12}$ are
uniform in $t_0$ and $t$,
these prove equicontinuity in the variable $z\in[0,1]$
of functions $\tilde{h}_{t_0,t}(w,z)$ in $H$.

In a similar way as the proof of equicontinuity with respect to $z$,
we have
\[
|\tilde{h}_{t_0,t}(w',z)-\tilde{h}_{t_0,t}(w,z)|\le 
I_{21}(w,w')+I_{22}(w,w'),
\]
where
\[
\arrb{l}\dsp
I_{21}(w,w')=C_h
\sum_{k\ge0} \int_{0=:u_k< u_{k-1}<u_{k-2}<\cdots<u_1<u_0\le t_0}
\\ \dsp \phantom{I_{21}(w,w')=} \times
|e^{-X(w')}-e^{-X(w)}|\,
\biggl(\prod_{i=0}^{k-1} (w'(1,u_i)+C_W)du_i\biggr)
du_{k-1}\,,
\arre
\]
with
\[
X(w)=\sum_{i=0}^{k-1} \Omega_{\theta,w,z}(u_{i+1},u_i)
+\Omega_{\theta,w,z}(u_0,t),
\]
and
\[
\arrb{l}\dsp
I_{22}(w,w')=C_h
e^{-\tilde{\Omega}_w(0,t)+C_Wt}
\sum_{k\ge1} \int_{0=:u_k< u_{k-1}<u_{k-2}<\cdots<u_1<u_0\le t_0}
\\ \dsp \phantom{I_{22}(w,w')=} \times
\biggl|
\prod_{i=0}^{k-1} \tilde{w}_{\theta,w',z}(u_{i+1},u_i)
-\prod_{i=0}^{k-1} \tilde{w}_{\theta,w,z}(u_{i+1},u_i)
\biggr|
\prod_{i=0}^{k-1} du_i\,.
\arre
\]
Note that \eqnu{omegaOmega} with \eqnu{Wnorm} implies
\[
|\Omega_{\theta,w',z}(u,v)-\Omega_{\theta,w,z}(u,v)|
\le \wnorm{w'-w} (v-u),\ \ 0\le u\le T.
\]
Using this, \eqnu{equicontiprf5}, and \eqnu{equicontiprf3}
in $I_{21}(w,w')$, we further have
\[
\arrb{l}\dsp
I_{21}(w,w')
\le
C_h e^{-\tilde{\Omega}_w(0,t)+C_Wt}\,e^{\wnorm{w'-w}t}
\wnorm{w'-w}t \,e^{\tilde{\Omega}_{w'}(0,t_0)+C_Wt_0}
\\ \dsp \phantom{I_{21}(w,w')}
\le
C_h e^{2C_WT}\,e^{\wnorm{w'-w}T}\wnorm{w'-w}T.
\arre
\]
With a similar argument, we also have
\[
\arrb{l}\dsp
I_{22}(w,w')
\le
C_he^{-\tilde{\Omega}_w(0,t)+C_Wt}
\sum_{k\ge1}\sum_{j=0}^{k-1}
\int_{0=:u_k< u_{k-1}<u_{k-2}<\cdots<u_1<u_0\le t_0}
\\ \dsp \phantom{I_{22}(w,w')\le} \times
\biggl(
\prod_{i=0}^{j-1} \tilde{w}_{\theta,w',z}(u_{i+1},u_i)
\biggr)\,
|\tilde{w}_{\theta,w',z}(u_{j+1},u_j)
-\tilde{w}_{\theta,w,z}(u_{j+1},u_j)|
\\ \dsp \phantom{I_{22}(w,w')\le} \times
\biggl(
\prod_{i=j+1}^{k-1} \tilde{w}_{\theta,w,z}(u_{i+1},u_i)
\biggr)
\prod_{i=0}^{k-1} du_i
\\ \dsp \phantom{I_{22}(w,w')}
\le C_he^{-\tilde{\Omega}_w(0,t)+C_Wt}\wnorm{w'-w}
\int_0^{t_0}  e^{\tilde{\Omega}_w(0,v)+C_Wv}
e^{\tilde{\Omega}_{w'}(v,t_0)+C_W(t_0-v)} \,dv
\\ \dsp \phantom{I_{22}(w,w')}
\le C_he^{-\tilde{\Omega}_w(0,t)+C_Wt}\wnorm{w'-w}
\\ \dsp \phantom{I_{22}(w,w')\le} \times
\biggl(
e^{\tilde{\Omega}_{w}(0,t_0)}t_0+
\int_0^{t_0} e^{\tilde{\Omega}_{w'}(0,v)}
\, |e^{\tilde{\Omega}_{w'}(v,t_0)}-e^{\tilde{\Omega}_{w}(v,t_0)}|
\,dv
\biggr)
\\ \dsp \phantom{I_{22}(w,w')}
\le
C_he^{2C_WT}\wnorm{w'-w}T
\,(1+ e^{\wnorm{w'-w}T}\wnorm{w'-w}T).
\arre
\]
Since the right hand sides of the bounds for $I_{21}$ and $I_{22}$ are
uniform in $t_0$ and $t$,
these prove equicontinuity of $H$ in $w\in W$.

This completes a proof of equicontinuity of $H$,
hence a proof of \thmu{fSRPLLN}.
\QED
\prfofe

\section{Hierarchy of multi time Gronwall inequality.}
\seca{mGronwall}

The following is a simple form of Gronwall's inequality.
\thmb
\thma{integopNr}
Let $T$ be a positive constant, and $a$ and $c$ be non-negative constants.
If $x:\ [0,T]\to\reals$ is an integrable function,
satisfying
\[
x(t)\le a+c\int_0^t x(s)\,ds, \ t\in[0,T],
\]
then
\begin{equation}
\eqna{integopNp2}
x(t)\le a\,e^{ct},\ t\in[0,T],
\end{equation}
holds.
\DDD\thme
The following is a generalization of \thmu{integopNr} to 
functions of more than $1$ variables,
where the case $q=1$ is \thmu{integopNr}.
\thmb
\thma{integopNqsum}
Let $T$ be a positive constant, $q$ a positive integer,
and $a$ and $c$ non-negative constants.
If $x:\ [0,T]^q\to\reals$ is an integrable function of $q$ variables,
satisfying
\[
\arrb{l}\dsp
x(t_1,\ldots,t_q)
\le a\,e^{c\,(t_1+\cdots+t_q)}\frac1q\sum_{i=1}^q e^{-ct_i}
\\ \dsp \phantom{x(t_1,\ldots,t_q)\le}
+\frac{c}q\sum_{i=1}^q \int_0^{t_i} (x(t_1,\ldots,t_q)|_{t_i=u})\,du,
\ \ (t_1,\ldots,t_q)\in [0,T]^q,
\arre
\]
then
\begin{equation}
\eqna{integopNqsum}
x(t_1,\ldots,t_q)\le a\,e^{c\,(t_1+\cdots+t_q)},
\ \ (t_1,\ldots,t_q)\in [0,T]^q,
\end{equation}
holds.
\DDD\thme
To prove \thmu{integopNqsum}, we start with the homogeneous case.
\thmb
\thma{integopqvar}
Let $T$ be a positive constant, $q$ a positive integer,
and $c$  a non-negative constant.
If $x:\ [0,T]^q\to\reals$ is an integrable function of $q$ variables,
satisfying
\begin{equation}
\eqna{integopqvarassump}
x(t_1,\ldots,t_q)\le c\sum_{i=1}^q \int_0^{t_i}
 (x(t_1,\ldots,t_q)|_{t_i=s})\,ds, \ \ (t_1,\ldots,t_q)\in [0,T]^q,
\end{equation}
then
\begin{equation}
\eqna{integopqvar}
x(t_1,\ldots,t_q)\le 0,\ \ (t_1,\ldots,t_q)\in [0,T]^q,
\end{equation}
holds.
\DDD
\thme
To prove \thmu{integopqvar}, we introduce a notation
\begin{equation}
\eqna{integopqvarmulti}
\arrb{l}\dsp
(A_{i,k} y)(t_1,\ldots,t_q)
\\ \dsp {}
=\left\{\arrb{ll}\dsp
\frac1{(k-1)!}\int_0^{t_i} (t_i-s)^{k-1}\, (y(t_1,\ldots,t_q))|_{t_i=s}\,ds,
& k=1,2,3,\ldots,\\ \dsp y(t_1,\ldots,t_q)
\ \ \ \ (\mbox{i.e., } A_{i,0}=\rm{id}), & k=0,
 \arre\right.
\arre
\end{equation}
for integrable function
$\dsp y:\ [0,T]^q\to\reals$ in $q$ variables and $i=1,\ldots,q$.
$A_{i,k}$, $k\in\pintegers$, $i=1,2,\ldots,q$, are commutative 
operators on the set of integrable functions.
In fact,
commutativity is obvious for $k=0$,
and by induction in $k$ we have
\begin{equation}
\eqna{integop}
A_{i,k} A_{i,\ell}=A_{i,1}^{k+\ell}=A_{i,k+\ell}=A_{i,\ell} A_{i,k},
\end{equation}
and Fubini's theorem implies for $k\ell>0$ and $i\ne j$ 
\[ \arrb{l}\dsp
(A_{i,k}A_{j,\ell}y)(t_1,\ldots,t_q)
\\ \dsp \phantom{}
=\frac{1}{(k-1)!}\frac{1}{(\ell-1)!}\int_0^{t_i}ds\,\int_0^{t_j}du\,
 (t_i-s)^{k-1}\, 
\\ \dsp \phantom{
 =\frac{1}{(k-1)!}\frac{1}{(\ell-1)!}\int_0^{t_i}ds\,\int_0^{t_j}du\,} \times
\left( (t_j-u)^{\ell-1}\, (y(t_1,\ldots,t_q))|_{t_j=u}
\right)|_{t_i=s}
\\ \dsp \phantom{}
=(A_{j,\ell}A_{i,k}y)(t_1,\ldots,t_q),
\arre \]
which prove
\begin{equation}
\eqna{integopcommutative}
A_{i,k}A_{j,\ell}=A_{j,\ell}A_{i,k},
\ \ k,\ell\in\pintegers,\ i,j\in\{1,\ldots,q\}.
\end{equation}
\lemb
\lema{integopqvar}
Under the assumptions of \thmu{integopqvar},
\begin{equation}
\eqna{integopqvarestim}
x(t_1,\ldots,t_q)\le c^N
 \sum_{(k_1,\ldots,k_q)\in\pintegers^q;\atop k_1+\cdots+k_q=N}
(A_{q,k_q}\,A_{q-1,k_{q-1}}\,\cdots\,A_{1,k_1}\, x)(t_1,\ldots,t_q),
\ N\in\pintegers,
\end{equation}
holds.
\DDD\leme
\prfb
The case $N=1$ of \eqnu{integopqvarestim}
is the assumption \eqnu{integopqvarassump} itself.
Assume that \eqnu{integopqvarestim} holds for some $N$.
Substituting \eqnu{integopqvarassump} in \eqnu{integopqvarestim},
and 
noting that sums, integrations, and multiplication of non-negative reals
have monotonicity,
we have
\[
x(t_1,\ldots,t_q)\le c^{N+1}
\sum_{i=1}^q 
 \sum_{(k_1,\ldots,k_q)\in\pintegers^q;\atop k_1+\cdots+k_q=N}
(A_{q,k_q}\,A_{q-1,k_{q-1}}\,\cdots\,A_{1,k_1}\,A_{i,1} x)(t_1,\ldots,t_q).
\]
Using \eqnu{integop} in the form $\dsp A_{i,k_i}A_{i,1} = A_{i,k_i+1}$,
we have \eqnu{integopqvar} for $N$ replaced by $N+1$.
\QED\prfe
\prfofb{\protect\thmu{integopqvar}}
For notational simplicity,
put $\vec{t}=(t_1,\ldots,t_q)$ in this proof.
The operator $A_{i,k}$ in \eqnu{integopqvarmulti} satisfies
\[
(A_{i,k} y)(\vec{t})\le \frac{t_i^k}{k!}\,
\sup_{\vec{t}\in[0,T]^q} y(\vec{t}),
\ \ \vec{t}\in [0,T]^q,
\]
for a integrable function $y$,
hence \eqnu{integopqvarestim} implies
\begin{equation}
\eqna{integopqvarprf}
x(\vec{t})\le
c^N \sum_{(k_1,\ldots,k_q)\in\pintegers^q;\atop k_1+\cdots+k_q=N}
\prod_{i=1}^q \frac{t_i^{k_i}}{k_i!}\,
\sup_{\vec{t}\in[0,T]^q} x(\vec{t}),
\ \ \vec{t}\in [0,T]^q,\ N\in\nintegers.
\end{equation}
For an arbitrary $\eps>0$, let
$\dsp\vec{t}_0=(t_{0,1},\ldots,t_{0,q})\in [0,T]^q$
be a vector (independent of $N$) such that
$\dsp x(\vec{t}_0)
\ge \sup_{\vec{t}\in [0,T]^q} x(\vec{t})-\eps$ holds.
Put
\[
a_N= c^N\sum_{(k_1,\ldots,k_q)\in\pintegers^q;\atop k_1+\cdots+k_q=N}
\prod_{i=1}^q \frac{t_{0,i}^{k_i}}{k_i!}.
\]
Then \eqnu{integopqvarprf} implies
\[
\sup_{\vec{t}\in [0,T]^q} x(\vec{t})
\le  x(\vec{t}_0)+\eps
\le a_N\sup_{\vec{t}\in[0,T]^q} x(\vec{t}) +\eps,
\]
hence
\[
\sup_{\vec{t}\in [0,T]^q} x(\vec{t})\,(1-a_N)\le \eps
\]
holds.
We see
\[
 \sum_{N=0}^{\infty} a_N
\le 
 \sum_{N=0}^{\infty} 
\sum_{(k_1,\ldots,k_q)\in\pintegers^q;\atop k_1+\cdots+k_q=N}
\prod_{i=1}^q \frac{(ct_{0,i})^{k_i}}{k_i!}
=\prod_{i=1}^q e^{ct_{0,i}}<\infty,
\]
so that, in particular, $\dsp \limf{N} a_N=0$,
which implies $\dsp 1-a_N\le \frac12$ for large $N$.
Hence
\[
\sup_{\vec{t}\in [0,T]^q} x(\vec{t})\le 2\eps,
\]
which proves \eqnu{integopqvar}.
\QED\prfofe

\prfofb{\protect\thmu{integopNqsum}}
Note that
\begin{equation}
\eqna{integopNqsum1}
x_1(s_1,\ldots,s_q)=a\,e^{c\,(s_1+\cdots+s_q)}
\end{equation}
satisfies
\begin{equation}
\eqna{integopNqsum2}
\arrb{l}\dsp
x_1(s_1,\ldots,s_q)
= e^{-cs_i} x_1(s_1,\ldots,s_q)
+c \int_0^{s_i} (x_1(s_1,\ldots,s_q)|_{s_i=u})\,du,
\\ \dsp
\ \ (s_1,\ldots,s_q)\in [0,t]^q,\ i=1,2,\ldots,q.
\arre
\end{equation}
Subtracting $x_1(s_1,\ldots,s_q)$ from \eqnu{integopNqsum},
and then using \eqnu{integopNqsum2}, we have
\[\arrb{l}\dsp
x(t_1,\ldots,t_q)-x_1(t_1,\ldots,t_q)
\le \frac{c}q\sum_{i=1}^q \int_0^{t_i}
 (x(t_1,\ldots,t_q)-x_1(t_1,\ldots,t_q))|_{t_i=s}\,ds,
\\ \dsp
 (t_1,\ldots,t_q)\in [0,T]^q,
\arre\]
which, with \thmu{integopqvar} and \eqnu{integopNqsum1}
implies \eqnu{integopNqsum}.
\QED
\prfofe

Finally, we give a result to be used in the proof of the main theorem
in \secu{MainProof} which contains recursion with respect to the
number of variables $q$ and a nonlinear term.
\thmb
\thma{integopNqrecs}
Let $T$ be a positive constant, $d$ be a non-negative constant
satisfying $d\le 1$,
and for each positive integer $q$
let $a_q$, $b_q$, and $c_q$ be non-negative constants.
Assume that, for a series of non-negative valued integrable functions
$\dsp x_q:\ [0,T]^q\to \preals$, $q\in\pintegers$, 
\begin{equation}
\eqna{integopNqrecsassump}
\arrb{l}\dsp
x_0=1,
\\ \dsp
x_q(t_1,\ldots,t_q)\le 
a_q\sum_{i=1}^q x_{q-1}(t_1,\ldots,\not{\!\!t_{i}},\ldots,t_q)^d
\\ \dsp \phantom{x_q(t_1,\ldots,t_q)\le }
+b_q\sum_{i=1}^q x_{q-1}(t_1,\ldots,\not{\!\!t_{i}},\ldots,t_q)
\\ \dsp \phantom{x_q(t_1,\ldots,t_q)\le }
+c_q\sum_{i=1}^q \int_0^{t_i} (x_q(t_1,\ldots,t_q)|_{t_i=s})\,ds,
\\ \dsp
(t_1,\ldots,t_q)\in [0,T]^q,\ \ q\in\nintegers,
\arre
\end{equation}
hold.
Put
\[
\tilde{c}_q=\max_{1\le k\le q} kc_k\,,\ q\in\nintegers,
\]
and define a sequence of non-negative constants $g_q$, $q=0,1,2,\ldots$,
recursively by
\[
g_0=1,\ \ g_q=q\,(a_q\,g_{q-1}^d+b_q\,g_{q-1}),\ q\in\nintegers.
\]
Then
\begin{equation}
\eqna{integopNqrecs}
x_q(t_1,\ldots,t_q)\le g_q\,e^{\tilde{c}_q\,(t_1+\cdots+t_q)},
\ \ 
(t_1,\ldots,t_q)\in [0,T]^q,\ \ q\in\nintegers,
\end{equation}
holds.
\DDD
\thme
\prfb
If $q=1$, \eqnu{integopNqrecsassump} reads
$\dsp
x_1(t_1)\le (a_1+b_1)+c_1\int_0^{t_1}x_1(s)\,ds,
$
hence \thmu{integopNr} implies
$\dsp x_1(t_1)\le g_1e^{c_1t_1}$, which proves
\eqnu{integopNqrecs} for $q=1$.

Let $q\ge 2$ and
assume that \eqnu{integopNqrecs} holds for $x_{q-1}$\,, as
\[
x_{q-1}(t_1,\ldots,t_{q-1})
\le g_{q-1} e^{\tilde{c}_{q-1}(t_1+\cdots+t_{q-1})}.
\]
This and \eqnu{integopNqrecsassump} for $x_q$ and $d\le 1$ imply
\[
\arrb{l}\dsp
x_q(s_1,\ldots,s_q)\le 
q(a_qg_{q-1}^d+b_qg_{q-1})\,e^{\tilde{c}_{q-1}\,(s_1+\cdots+s_q)}
\frac1q \sum_{i=1}^q e^{-\tilde{c}_{q-1}s_i}
\\ \dsp \phantom{x_q(s_1,\ldots,s_{q})\le }
+qc_q\frac1q\sum_{i=1}^q \int_0^{s_i} (x_q(s_1,\ldots,s_q)|_{s_i=u})\,du,
\arre
\]
which, with \thmu{integopNqsum}, implies \eqnu{integopNqrecs} for $x_q$.
\QED
\prfe

\section{Proof of the main theorem.}
\seca{MainProof}

\subsection{Convergence of the spatial distribution function.}
Here we will prove the essential part of the infinite particle limit,
the convergence of spatial distribution function.

In analogy to \eqnu{JtopnuKtheta},
define,
for each $i=1,2,\ldots,N$ and $0\le t_0\le t\le T$ and $0\le y_0\le 1$,
\begin{equation}
\eqna{JtopnuK}
J\ssN_i(t_0,t)=\{ \omega\in\Omega\mid 
\tilde{\nu}\ssN_i(t)(\omega)>\tilde{\nu}\ssN_i(t_0)(\omega)\}.
\end{equation}
By similar arguments as for \eqnu{varphiNthetaJ} and \eqnu{yCNtheta},
$\varphi\ssN$ in \eqnu{varphiN} and $Y\ssN_C$ in \eqnu{varphiNYNC}
respectively satisfies
\begin{equation}
\eqna{varphiNJ}
\varphi\ssN(dw,(y_0,t_0),t)=
\frac1N\sum_{j;\ Y\ssN_j(t_0)\ge y_0}
 \chrfcn{J\ssN_j(t_0,t)^c}\,\delta_{w_j}(dw),
\end{equation}
and
\begin{equation}
\eqna{yCN}
Y\ssN_C((y_0,t_0),t)
=y_0+\frac1N\sum_{j;\ Y\ssN_j(t_0)\ge y_0}
 \chrfcn{J\ssN_j(t_0,t)}\,.
\end{equation}
%
\thmb
\thma{SRPLLN}
Assume \eqnu{Cw}, \eqnu{muN0conv}, \eqnu{rw}, and \eqnu{lambda1}.
Then there exists $\delta'>0$ and an integer $p_0$ 
satisfying $2p_0\delta'>1$, such that for any integer $p\ge p_0$
there exists a positive constant $C$
depending only on $p$ and $\delta'$,
(independent of $N$ and $h$,) such that
for any bounded continuous $h:\ W\to\reals$
\begin{equation}
\eqna{SRPHDLepsilon}
\EE{ \sup_{(\gamma,t)\in\Delta_T}\biggl|
\varphi\ssN(h,\gamma,t)-\varphi_{y_C}(h,\gamma,t)
\biggr|^{2p} }
\le \frac{C\,C_h^{2p}}{N^{2p\delta'}}\,,
\ N\in\nintegers,
\end{equation}
holds, where $C_h$ is as in \eqnu{hdata}.
\DDD
\thme
Note that \eqnu{varphiNYNC} and \eqnu{yCFPG} with $\theta=y_C$ imply
\[
Y\ssN_C(\gamma,t)-y_C(\gamma,t)=\frac{[N(1-y_0)]}N-(1-y_0)
+\varphi_{y_C}(W,\gamma,t)-\varphi\ssN(W,\gamma,t).
\]
Applying $\dsp (a+b)^{2p}\le2^{2p-1}(a^{2p}+b^{2p})$,
valid for $a,b\ge0$ and $2p\ge1$,
\thmu{SRPLLN} with $h(w)=1$, $w\in W$, therefore implies
\begin{equation}
\eqna{YCNyCepsilon}
\EE{ \sup_{(\gamma,t)\in\Delta_T}\biggl|
Y\ssN_C(\gamma,t)-y_C(\gamma,t)
\biggr|^{2p} }
\le \frac{2^{2p-1}C}{N^{2p\delta'}}+\frac{2^{2p-1}}{N^{2p}}\,,
\ N\in\nintegers,
\end{equation}
with the assumptions and notations of the Theorem.

\subsubsection{Coupling of the original and the flow driven model.}

In view of \thmu{fSRPLLN}, it suffices to prove the following
for \thmu{SRPLLN} to hold.
\thmb
\thma{SRPLLNdiffmoment}
Assume \eqnu{Cw}, \eqnu{muN0conv}, \eqnu{rw}, and \eqnu{lambda1},
and let $\delta$ be as in \eqnu{muN0conv}.
Then there exists $\delta'>0$ and an integer $p_0$ 
satisfying $2p_0\delta'>1$, such that for any integer $p\ge p_0$
there exists a positive constant $C$, (independent of $N$ and $h$,)
such that for any bounded continuous $h:\ W\to\reals$
\begin{equation}
\eqna{SRPLLNdiffmoment}
\EE{ \sup_{((y_0,t_0),t)\in\Delta_T}\biggl|
\varphi\ssN(h,\gamma,t)-\varphi\ssNyC(h,\gamma,t)
\biggr|^{2p} }
\le \frac{C\,C_h^{2p}}{N^{2p\delta'}}\,,
\ N\in\nintegers,
\end{equation}
where $C_h$ is as in \eqnu{hdata}.
\DDD
\thme

For $t\in[0,T]$ and $i\in\{1,2,\ldots,N\}$, put
\begin{equation}
\eqna{omegawedge}
\tilde{w}\ssN_{i,\wedge}(t)=
 w_i(Y\ssN_i(t-),t)\wedge w_i(y_C(\gamma\ssNyC_i(t-),t),t),
\end{equation}
and
\begin{equation}
\eqna{omegavee}
\tilde{w}\ssN_{i,\vee}(t)=
 w_i(Y\ssN_i(t-),t)\vee w_i(y_C(\gamma\ssNyC_i(t-),t),t),
\end{equation}
and denote the event that the $i$-th particle 
$Y\ssN_i(s)$ of \eqnu{SRPY} and $Y\ssNyC_i(s)$ of \eqnu{SRPYtheta}
jump to top at same times in the interval $(t_0,t]$ by
\begin{equation}
\eqna{MRUjumpseqdiffnu}
\mathcal{K}\ssN_i(t_0,t)=\{ \omega\in\Omega\mid 
\nu\ssN_i(\{(s,\xi)\mid
\tilde{w}\ssN_{i,\wedge}(s)<\xi\le \tilde{w}\ssN_{i,\vee}(s),
\ s\in(t_0,t]\ \}) =0  \}.
\end{equation}

Fix a bounded continuous function $h:\ W\to\reals$,
and let $C_h$ be as in \eqnu{hdata}.
Using the definitions \eqnu{varphiNJ} and \eqnu{varphiNthetahexplicit},
put
\begin{equation}
\eqna{diffmomentb}
\arrb{l}\dsp
\Delta\varphi\ssN(\gamma,t)
=\varphi\ssN(h,\gamma,t)-\varphi\ssNtheta(h,\gamma,t)
\\ \dsp \phantom{\Delta\varphi\ssN(\gamma,t)}
=\frac1N\sum_{j;\ y\ssN_j\ge y_0} h(w_j) 
(\chrfcn{J\ssN_j(t_0,t)^c}- \chrfcn{J\ssNtheta_j(t_0,t)^c}),
\\ \dsp {}
\gamma=(y_0,t_0)\in\Gamma_t,\ t\in[0,T].
\arre
\end{equation}
Then 
\eqnu{MRUjumpseqdiffnu} and \eqnu{hdata} imply,
\[
|\Delta\varphi\ssN(\gamma,t)|
\le \frac{C_h}N\sum_{j=1}^N \chrfcn{\mathcal{K}\ssN_j(t_0,t)^c},
\ \ \gamma=(y_0,t_0)\in\Gamma_t,\ t\in[0,T].
\]
The monotonicity of $\mathcal{K}\ssN_j(t_0,t)^c$ with respect
to $t$ and $t_0$ further implies
\begin{equation}
\eqna{diffmomentnu}
\sup_{(\gamma,t)\in\Delta_T} |\Delta\varphi\ssN(\gamma,t)|
\le
\frac{C_h}N
\sum_{i=1}^N \chrfcn{\mathcal{K}\ssN_{i}(0,T){}^c},
\end{equation}
Proof of \thmu{SRPLLNdiffmoment} therefore reduces to evaluation of 
the event $\mathcal{K}\ssN_{i}(0,T)$.

\subsubsection{Event with different jumps to top.}
\seca{eventcmp}

As an analog of \eqnu{latestarrivaltimeSItilde},
define a sequence of stopping times,
$0=\tau\ssN_{i,0}<\tau\ssN_{i,1}<\cdots$,
by
\begin{equation}
\eqna{coupledjumptimeiSRP}
\arrb{l}\dsp
\tau\ssN_{i,0}=0,
\\ \dsp
\tau\ssN_{i,k+1}=\inf\{t>\tau\ssN_{i,k}\mid 
\nu\ssN_i(
\{(s,\xi)\in(\tau\ssN_{i,k},T]\times\preals\mid
\\ \dsp \phantom{\tau\ssN_{i,k+1}=}
0\le \xi\le w_i(Y\ssN_i(s-),s)\})>0\},\ k\in\pintegers.
\arre
\end{equation}
$\tau\ssN_{i,k}$ is the time that the particle $i$ in the
(original) stochastic ranking process jumps to the top for 
the $k$-th time.
A corresponding analog of \eqnu{flowstarttime} is
 \begin{equation}
\eqna{flowstarttimeSRP}
\arrb{l}\dsp
\gamma\ssN_i(t)=\left\{\arrb{ll}\dsp 
(y\ssN_i,0),& 0\le t<\tau\ssN_{i,1},
\\ \dsp 
(0,\tau\ssN_{i,k}),& \tau\ssN_{i,k}\le t<\tau\ssN_{i,k+1},
\ k=1,2,\ldots.
\arre\right.
\arre
 \end{equation}
A property corresponding to \eqnu{YNthetaieqYNCfg} then is
\begin{equation}
\eqna{YNieqYNCfg}
Y\ssN_i(t)=Y\ssN_C(\gamma\ssN_i(t),t),\ t\in[0,T],
\end{equation}
which can be proved in a similar way as a proof of \eqnu{YNthetaieqYNCfg}
in \lemu{YNCeqYNi}.
This decomposition in particular decomposes the dependence as random
variables; if we temporarily denote by $X\in \mathcal{F}$,
a fact that a random variable $X:\ \Omega\to\reals$ 
is $\mathcal{F}$-measurable, and denote by $\sigma[Z]$ a sigma algebra
generated by a random variable $Z$, we have
\begin{equation}
\eqna{YNisigma}
\arrb{l}\dsp
Y\ssN_C((y_0,t_0),t)\in\sigma[\{\nu\ssN_j\mid Y\ssN_j(t_0)>y_0\}],
\\ \dsp
\gamma\ssN_i\in \sigma[\{\tau\ssN_{i,k}\wedge t\mid k\in\nintegers\}].
\arre
\end{equation}

Define an analog of the stopping times \eqnu{coupledjumptimeiSRP} 
using \eqnu{omegawedge} by
\begin{equation}
\eqna{coupledjumptimeiwedgeSRP}
\arrb{l}\dsp
\tau\ssN_{i,\wedge,0}=0,
\\ \dsp {}
\tau\ssN_{i,\wedge,k}=\inf\{t>\tau\ssN_{i,\wedge,k-1}\mid
\nu\ssN_i(\{(\xi,s)\mid 0\le \xi\le \tilde{w}\ssN_{i,\wedge}(s),
\ 0\le s\le t\})>0\},\ 
\\ \dsp {}
k\in\nintegers,
\arre
\end{equation}
and denote by $\sigma\ssN_i$, the time that the particle pair with label 
$i$ of the original model and the flow driven model have different jumps
to the top for the first time;
\begin{equation}
\eqna{difftime}
\sigma\ssN_i=\inf\{t\in[0,T]\mid \mathcal{K}\ssN_i(0,t)^c\}.
\end{equation}
The definition implies
\begin{equation}
\eqna{wedgeeqoriginal}
\tau\ssN_{i,k}<\sigma\ssN_i\ \Rightarrow
\ \tau\ssN_{i,\wedge,k}=\tau\ssN_{i,k}=\tau\ssNyC_{i,k}\,,
\end{equation}
where $\tau\ssNyC_{i,k}$ is defined in \eqnu{latestarrivaltimeSItilde},
with $\theta=y_C$\,.
Indepndence of $\nu\ssN_i(A)$ and $\nu\ssN_i(B)$ for
the exclusive events $A$ and $B$ implies
\begin{equation}
\eqna{wedge_difference_indep}
\{\tau\ssN_{i,\wedge,k}\mid k\in\pintegers\}\ \perp\ \sigma\ssN_i\,.
\end{equation}

Using \eqnu{omegavee}, \eqnu{omegawedge}, \eqnu{Cw}C
and $a\vee b=|a-b|+a\wedge b$,
we have
\begin{equation}
\eqna{recintdec}
\arrb{l}\dsp
\{(\xi,s)\in\nreals2 \mid
\tilde{w}\ssN_{i,\wedge}(s)<\xi \le  \tilde{w}\ssN_{i,\vee}(s),\ 0< s\le t\}
\\ \dsp \phantom{}
\subset\{(\xi,s)\in\nreals2 \mid
0\le \xi-\tilde{w}\ssN_{i,\wedge}(s) \le 
 C_W |Y\ssN_{i }(s)- y_C(\gamma\ssNyC_{i }(s),s)|,
\\ \dsp \phantom{\subset\{(\xi,s)\in\nreals2 \mid}
\ 0< s\le t\}
\\ \dsp \phantom{}
\subset
\bigcup_{k=1}^{\infty} \{(\xi,s)\in\nreals2 \mid
0\le \xi-\tilde{w}\ssN_{i,\wedge}(s) \le 
 C_W |Y\ssN_{i }(s)-y_C(\gamma\ssNyC_{i }(s),s)|,
\\ \dsp \phantom{\subset\bigcup_{k=1}^{\infty}\{(\xi,s)\in\nreals2 \mid}
\ t\wedge\tau\ssN_{i ,\wedge,k-1}< s\le t\wedge\tau\ssN_{i ,\wedge,k}\}.
\arre
\end{equation}
Note that for each $i$
\[
Y\ssN_i(\tau\ssN_{i,k-1})=Y\ssNyC_i(\tau\ssN_{i,k-1})=0,\ \mbox{ on }
\tau\ssN_{i,k-1}<\sigma\ssN_i\,.
\]
Note also that the definition \eqnu{difftime} implies
$\dsp\{ \sigma\ssN_i>s\} =\mathcal{K}\ssN_i(0,s)$.
Hence, \eqnu{YNieqYNCfg}, \eqnu{YNthetaieqYNCfg}, \eqnu{yCN}, and
\eqnu{yCNtheta} imply, 
with similar arguments for deriving 
\eqnu{diffmomentnu} from \eqnu{MRUjumpseqdiffnu},
\[ \arrb{l}\dsp
|Y\ssN_i(s)-Y\ssNyC_i(s)|
=|Y\ssN_C(\gamma\ssN_i(s),s)-Y\ssNyC_C(\gamma\ssN_i(s),s)|
\\ \dsp \phantom{|Y\ssN_i(s)-Y\ssNyC_i(s)|}
=\frac1N|\sum_{j\ne i}(  \chrfcn{(J\ssN_j(\tau\ssN_{i,k-1},s))^c}
-\chrfcn{(J\ssN_j(\tau\ssNyC_{i,k-1},s))^c}  )|
\\ \dsp \phantom{|Y\ssN_i(s)-Y\ssNyC_i(s)|}
\le
\frac1N\sum_{j\ne i} \chrfcn{\mathcal{K}\ssN_j(\tau\ssN_{i,k-1},s)^c}
\\ \dsp
\mbox{ on }
\mathcal{K}\ssN_{i}(0,s),
\ \tau\ssN_{i,k-1}\le s<\tau\ssN_{i,k}.
\arre \]
This with \eqnu{recintdec} then implies
\begin{equation}
\eqna{Ydiffctbd}
\arrb{l}\dsp
\mathcal{K}\ssN_i(0,t)^c
\\ \dsp {}
\subset\bigcup_{k=1}^{\infty}\biggl(
\{ \omega\in\Omega\mid \nu\ssN_i(\{(s,\xi)\mid
0\le \xi-\tilde{w}\ssN_{i,\wedge}(s)
\\ \dsp \phantom{\subset\bigcup_{k=1}^{\infty}\biggl(
\{ \omega\in\Omega\mid \nu\ssN_i(\{(s,\xi)\mid }
\le 
 C_W |Y\ssNyC_{i }(s)-y_C(\gamma\ssNyC_{i }(s),s)|
\\ \dsp \phantom{\subset\bigcup_{k=1}^{\infty}\biggl(
\{ \omega\in\Omega\mid \nu\ssN_i(\{(s,\xi)\mid \le \ }
+\frac{C_W}N\sum_{j\ne i}\chrfcn{\mathcal{K}\ssN_j(\tau\ssN_{i,k-1},s)^c},
\\ \dsp \phantom{\subset\bigcup_{k=1}^{\infty}\biggl(
\{ \omega\in\Omega\mid \nu\ssN_i(\{(s,\xi)\mid \ }
\ t\wedge\tau\ssN_{i ,\wedge,k-1}< s\le t\wedge\tau\ssN_{i ,\wedge,k}
\})>0\ \}
\\ \dsp \phantom{\subset\bigcup_{k=1}^{\infty}\biggl(}
\cap \mathcal{K}\ssN_i(0,\tau\ssN_{i,k-1})\biggr).
\arre
\end{equation}

\subsubsection{Application of Gronwall hierarchy.}

For 
$q=1,2,\ldots,N$ and $t_i\in[0,T]$, $i=1,\ldots,q$, put
\begin{equation}
\eqna{Xq}
X\ssN_q(t_1,\ldots,t_q)=\max_{\{i_1,\ldots,i_q\}\subset\{1,\ldots,N\}}
\EE{\prod_{\alpha=1}^q \chrfcn{\mathcal{K}\ssN_{i_{\alpha}}(0,t_{\alpha})^c}}.
\end{equation}

\eqnu{diffmomentb}, \eqnu{diffmomentnu}, and \eqnu{Xq} imply that
to prove \thmu{SRPLLNdiffmoment}, it suffices to find
$\delta'>0$ and integer $p_0$ satisfying $2p_0\delta'>1$,
such that for any integer $p\ge p_0$,
\begin{equation}
\eqna{SRPLLNdiffmomentbsuff}
\frac1{N^{2p}}\sum_{q=1}^{2p}\Ccomb{N}{q} d(2p,q) X\ssN_q(T,\ldots,T)
\le \frac{C}{N^{2p\delta'}},
\end{equation}
for some $C>0$ independent of $N$.
Here, $d(r,q)$ is the number of surjections from a finite set of size $r$
to a set of size $q$, which is determined inductively by
\begin{equation}
\eqna{no_surj}
d(r,1)=1,
\ \mbox{and}\ d(r,q)=q^r-\sum_{k=1}^{q-1} \Ccomb{q}k d(r,q-k),
\ q=2,3,\ldots,r.
\end{equation}

Fix $q$ and $\{i_1,\ldots,i_q\}$ in the right hand side of \eqnu{Xq}.
Let $\alpha\in\{1,\ldots,q\}$ be the suffix such that $y\ssN_{i_\alpha}$
is the smallest among $y\ssN_{i_1}$, $\ldots$, $y\ssN_{i_q}$, and
put $i_0=i_{\alpha}$.
At times $\tau\ssN_{i_0,k}$, $k\in\pintegers$, the particle $i_0$
is at the top position, namely, for $i_{\alpha}\ne i_0$,
\[ \arrb{l}\dsp
Y\ssN_{i_0}(\tau\ssN_{i_0,k})=0< Y\ssN_{i_{\alpha}}(\tau\ssN_{i_0,k}),
\ k=1,2,\ldots,
\\ \dsp 
Y\ssN_{i_0}(\tau\ssN_{i_0,0})=Y\ssN_{i_0}(0)=y\ssN_{i_0}
< y\ssN_{i_{\alpha}}=Y\ssN_{i_{\alpha}}(\tau\ssN_{i_0,0}).
\arre\]
Hence up to the first jump to the top,
each $Y\ssN_{i_{\alpha}}(t)$ with $i_{\alpha}\ne i_0$ is independent of
$\nu_{i_0}$.
Therefore,
\[ \arrb{l}\dsp
\EE{\prod_{\alpha=1}^q \chrfcn{\mathcal{K}\ssN_{i_{\alpha}}(0,t_{\alpha}){}^c}}
\\ \dsp {}
=\EE{\prod_{\alpha;\ i_{\alpha}\ne i_0}
 \chrfcn{(\mathcal{K}\ssN_{i_{\alpha}}(0,t_{\alpha}))^c}
 \prb{ (\mathcal{K}\ssN_{i_0}(0,t_0))^c
\mid \{\nu_j,\ j\ne i_0\}\cup \{\tau\ssN_{i_0,\wedge,k}\}}},
\arre \]
where
$\dsp\prb{\cdot \mid \{\nu_j,\ j\ne i_0\}\cup \{\tau\ssN_{i_0,\wedge,k}\}}$
denotes conditional probability conditioned on the sigma algebra
generated by 
$\nu_j$, $j\ne i_0$, and $\tau\ssN_{i_0,\wedge,k}$, $k\in\nintegers$.
This with \eqnu{wedge_difference_indep}, \eqnu{wedgeeqoriginal},
\eqnu{Ydiffctbd}, and 
\[
\prb{\nu(A)>0}=1-e^{-|A|}\le |A|
\]
for a unit Poisson random measure $\nu$, further leads to
\begin{equation}
\eqna{diffprodSIbdprobrec}
\arrb{l}\dsp
\EE{\prod_{\alpha=1}^q \chrfcn{\mathcal{K}\ssN_{i_{\alpha}}(0,t_{\alpha}){}^c}}
\\ \dsp {}
\le \EE{\prod_{\alpha;\ i_{\alpha}\ne i_0}
 \chrfcn{(\mathcal{K}\ssN_{i_{\alpha}}(0,t_{\alpha}))^c}
\sum_{k=1}^{\infty}
 \chrfcn{\mathcal{K}\ssN_{i_0}(0,t_0\wedge \tau\ssN_{i_0,k-1})}
\\ \dsp \phantom{=\EE{\prod_{\alpha;\ i_{\alpha}\ne i_0}} }\times
\biggl(
 \int_{t_0\wedge\tau\ssN_{i_0,\wedge,k-1}}^{%
t_0\wedge\tau\ssN_{i_0,\wedge,k}}
 C_W |Y\ssNyC_{i_0 }(s)-y_C(\gamma\ssNyC_{i_0 }(s),s)|
\,ds
\\ \dsp \phantom{=\EE{\prod_{\alpha;\ i_{\alpha}\ne i_0}} \biggl( \ }
+ \int_{t_0\wedge\tau\ssN_{i_0,\wedge,k-1}}^{%
t_0\wedge\tau\ssN_{i_0,\wedge,k}}
\frac{C_W}N\sum_{j\ne i_0}
\chrfcn{\mathcal{K}\ssN_j(t_0\wedge\tau\ssN_{i,k-1},s)^c}\,ds
\biggr) \ }
\\ \dsp {}
\le\frac{C_W}N\,\sum_{j=1}^N
 \EE{\prod_{\alpha;\ i_{\alpha}\ne i_0}
 \chrfcn{(\mathcal{K}\ssN_{i_{\alpha}}(0,t_{\alpha}))^c}
\sum_{k=1}^{\infty}
\chrfcn{\mathcal{K}\ssN_{i_0}(0,t_0\wedge\tau\ssN_{i_0,\wedge,k-1})}
\\ \dsp \phantom{=\EE{\prod_{\alpha;\ i_{\alpha}\ne i_0}} }\times
\int_{t_0\wedge\tau\ssN_{i_0,\wedge,k-1}}^{t_0\wedge\tau\ssN_{i_0,\wedge,k}} 
\chrfcn{\mathcal{K}\ssN_j(t_0\wedge\tau\ssN_{i_0,k-1},s-)^c}
\,ds}
\\ \dsp \phantom{\le}
+
C_W\EE{\prod_{\alpha;\ i_{\alpha}\ne i_0}
 \chrfcn{(\mathcal{K}\ssN_{i_{\alpha}}(0,t_{\alpha}))^c}
\sum_{k=1}^{\infty}
\chrfcn{\mathcal{K}\ssN_{i_0}(0,t_0\wedge\tau\ssN_{i_0,\wedge,k-1})}
\\ \dsp \phantom{=\EE{\prod_{\alpha;\ i_{\alpha}\ne i_0}} } \times
\int_{t_0\wedge\tau\ssN_{i_0,\wedge,k-1}}^{t_0\wedge\tau\ssN_{i_0,\wedge,k}} 
|Y\ssNyC_{i_0}(s)-y_C(\gamma\ssNyC_{i_0}(s),s)|\,ds}.
\arre
\end{equation}

Using
$\dsp\chrfcn{\mathcal{K}\ssN_{i_0}(0,t_0\wedge\tau\ssN_{i_0,k-1})}\le1$ and
$\dsp\chrfcn{\mathcal{K}\ssN_j(t_0\wedge\tau\ssN_{i_0,k-1},s-)^c}\le
\chrfcn{\mathcal{K}\ssN_j(0,s-)^c}$,
we see that for $i,j=1,2,\ldots,N$,
\begin{equation}
\eqna{diffdecomp1contribd}
\sum_{k\ge1}\chrfcn{\mathcal{K}\ssN_i(0,t\wedge\tau\ssN_{i,k-1})}
\int_{t\wedge\tau\ssN_{i,k-1}}^{t\wedge\tau\ssN_{i,k}} 
\chrfcn{\mathcal{K}\ssN_j(t\wedge\tau\ssN_{i,k-1},s-)^c}\,ds
\le \int_0^t \chrfcn{\mathcal{K}\ssN_j(0,s-)^c}\,ds.
\end{equation}
Substituting \eqnu{diffdecomp1contribd} in the first term of 
the right hand side of \eqnu{diffprodSIbdprobrec},
and bounding the characteristic function on the right hand side by $1$
for $j\in\{i_1,\ldots,i_q\}\setminus\{i_0\}$,
and bounding the characteristic function for $i_0$ also by $1$ 
in the second term of 
the right hand side of \eqnu{diffprodSIbdprobrec},
we have
\begin{equation}
\eqna{integopqvarprf2}
\arrb{l}\dsp
X\ssN_q(t_1,\ldots,t_q)
\\ \dsp{}
\le 
\frac{C_WT(q-1)}N \sum_{i_0=1}^q
 X\ssN_{q-1}(t_1,\ldots,\not{\!\!t_{i_0}},\ldots,t_q)
\\ \dsp{}
+C_W  \sum_{i_0=1}^q\int_0^{t_{i_0}} X\ssN_q(t_1,\ldots,s,\ldots,t_q)\,ds
\\ \dsp{}
+C_W
\max_{\{i_1,\ldots,i_q\}\subset\{1,\ldots,N\}}\int_0^{t_{i_0}}
\\ \dsp \phantom{{}+}
\EE{\prod_{i_{\alpha}\ne i_0}
 \chrfcn{\mathcal{K}_{i_{\alpha}}(0,t_{\alpha})^c}
\,|Y\ssNyC_{i_{i_0}}(s)-y_C((\gamma\ssNyC_{i_{i_0}}(s),s)|}\,ds.
\arre
\end{equation}
Here,
$\dsp X\ssN_{q-1}(t_1,\ldots,\not{\!\!t_{i_0}},\ldots,t_q)$ is
the function in \eqnu{Xq} with $q$ replaced by $q-1$ and with
$q-1$ variables obtained by excluding $t_{i_0}$ from $t_1,\ldots,t_q$,
and the variables for
$\dsp X\ssN_q(t_1,\ldots,s,\ldots,t_q)$ is $t_1,\ldots,t_q$
with $t_{i_0}$ replaced by $s$.
Applying H\"older's inequality in the form
\[
\EE{|X\,Y|}\le \EE{|X|^{2p/(2p-1)}}^{1-(2p)^{-1}}\EE{|Y|^{2p}}^{1/(2p)}
\]
to the last term in the right hand side of \eqnu{integopqvarprf2},
and using \eqnu{YNthetaieqYNCfg} and \eqnu{fSRPyCHDLepsilon}, we have
\[
\arrb{l}\dsp
X\ssN_q(t_1,\ldots,t_q)
\\ \dsp{}
\le 
\frac{C_WT(q-1)}N 
 \sum_{i_0=1}^q X\ssN_{q-1}(t_1,\ldots,\not{\!\!t_{i_0}},\ldots,t_q)
\\ \dsp\phantom{\le}
+C_W  \sum_{i_0=1}^q\int_0^{t_{i_0}} X\ssN_q(t_1,\ldots,s,\ldots,t_q)\,ds
\\ \dsp\phantom{\le}
+\frac{C_WTC}{N^{\delta}} \sum_{i_0=1}^q
 (X\ssN_{q-1}(t_1,\ldots,\not{\!\!t_{i_0}},\ldots,t_q))^{(2p-1)/(2p)},
\\ \dsp{}
q=1,2,\ldots,
\\ \dsp{}
X\ssN_0=1.
\arre
\]
Applying \thmu{integopNqrecs}, with
$a_q=C_WTCN^{-\delta}$, $b_q=C_WT(q-1)N^{-1}$, $c_q=C_W$, $d=1-(2p)^{-1}$,
we have
\begin{equation}
\eqna{Xqbd}
\arrb{l}\dsp
X\ssN_q(t_1,\ldots,t_q)
\le g_q e^{q^2C_WT},\ \ t_i\in[0,T],\ i=1,\ldots,q,
\ q\in\nintegers;
\\ \dsp
g_0=1,\ \ g_q=qC_WT(C N^{-\delta} g_{q-1}^d+(q-1)N^{-1}g_{q-1}),
\ q\in\nintegers.
\arre
\end{equation}
For large $N$ we have $g_{q-1}<1$, hence with $d<1$, 
we further have $g_{q-1}^d>g_{q-1}$, 
and also with $0<\delta<1$ for large $N$ we have $N^{-\delta}>(q-1)N^{-1}$,
so that
\[
g_q\le qC_WTg_{q-1}^dN^{-\delta}(C+1),\ q\in\nintegers.
\]
By induction in $q$,
\[
 g_q\le q!(C_WT\vee1)^q(C+1)^q\frac1{N^{2p\delta(1-(1-(2p)^{-1})^q)}}\,,
\ q\in\nintegers.
\]
Since $\dsp 1-(1-\frac1{2p})q$ is decreasing in $q$, we therefore have
\begin{equation}
\eqna{Xpbd}
\arrb{l}\dsp
X\ssN_{q}(t_1,\ldots,t_{q})\le
q!(C_WT\vee1)^q(C+1)^qe^{q^2C_WT}
\frac1{N^{2p\delta(1-(1-(2p)^{-1})^{2p})}}\,,\ 
\\ \dsp
q=1,2,\ldots,p.
\arre
\end{equation}

Choose $\delta'$ to be any positive constant satisfying
\[
0<\delta'<(1-\frac1e)\delta.
\]
Since $\dsp \limf{p} (1-\frac1{2p})^{2p}=e^{-1}<1$,
there exists an integer $\dsp p_0>\frac1{2\delta'}$ such that
\[
\delta'<(1-(1-\frac1{2p})^{2p})\delta ,\ \ p= p_0,\ p_0+1,\ \ldots.
\]
With \eqnu{Xpbd} we arrive at
\[
X\ssN_q(t_1,\ldots,t_q)\le
q!(C_WT\vee1)^q(C+1)^qe^{q^2C_WT}
\frac1{N^{2p\delta'}},\ q=1,2,\ldots,2p,
\]
for $p=p_0,p_0+1,\ldots$.
Since $\dsp \Ccomb{N}q \le \frac{N^q}{q!}\le \frac{N^{2p}}{q!}$\,,
this proves \eqnu{SRPLLNdiffmomentbsuff}.

This completes a proof of \thmu{SRPLLNdiffmoment},
and therefore, of \thmu{SRPLLN}.

\subsection{Proof of \protect\thmu{HDLposdep}.}

Let $y\in[0,1]$ and let
$h:\ W\to\reals$ be a bounded continuous function with $C_h>0$ as
in \eqnu{hdata}.
For the flow
$y_C\in\Theta_T$ in \thmu{FPtheorem}, the definition of $\Theta_T$
in \eqnu{class_of_conti_flows} implies that for each $t\in[0,T]$,
$\Gamma_t\ni\gamma\mapsto y_C(\gamma,t)\in[0,1]$ is surjective.
Therefore there exists $\gamma_t:\ [0,1]\to \Gamma_t$ such that
\begin{equation}
\eqna{HDLposdepprf}
y_C(\gamma_t(y),t)=y,\ y\in[0,1].
\end{equation}
We then have, using \eqnu{varphiyC}, \eqnu{varphithetah},
\eqnu{HDLposdepprf},
\eqnu{varphiN}, \eqnu{yCFPG}, \eqnu{hdata}, \eqnu{mut},
and \eqnu{HDLposdepprf} in turn, 
\[ \arrb{l}\dsp 
\biggl|
\int_W h(w) \mu\ssN_t(dw\times[y,1])-\int_W h(w) \mu_t(dw\times[y,1])
\biggr|
\\ \dsp {}
=
\biggl|
\int_W h(w) \mu\ssN_t(dw\times[y,1])-\varphi_{y_C}(h,\gamma_t(y),t)
\biggr|
\\ \dsp {}
\le 
\sup_{\gamma\in\Gamma_t} |\varphi\ssN(h,\gamma,t)-\varphi_{y_C}(h,\gamma,t)|
\\ \dsp \phantom{\le}
+\biggl|
\int_W h(w) \mu\ssN_t(dw\times[y,1])-
\int_W h(w) \mu\ssN_t(dw\times[Y\ssN_C(\gamma_t(y),t),1])
\biggr|
\\ \dsp {}
\le 
\sup_{\gamma\in\Gamma_t} |\varphi\ssN(h,\gamma,t)-\varphi_{y_C}(h,\gamma,t)|
+C_h |y-Y\ssN_C(\gamma_t(y),t)|
\\ \dsp \phantom{}
\le 
\sup_{\gamma\in\Gamma_t} |\varphi\ssN(h,\gamma,t)-\varphi_{y_C}(h,\gamma,t)|
+C_h |y_C(\gamma_t(y),t)-Y\ssN_C(\gamma_t(y),t)|
\\ \dsp \phantom{}
\le 
\sup_{\gamma\in\Gamma_t} |\varphi\ssN(h,\gamma,t)-\varphi_{y_C}(h,\gamma,t)|
+C_h \sup_{\gamma\in\Gamma_t} |y_C(\gamma,t)-Y\ssN_C(\gamma,t)|.
\arre \]
\thmu{SRPLLN} and \eqnu{YCNyCepsilon} then imply
\begin{equation}
\eqna{HDLposdepprf3}
\arrb{l}\dsp
\EE{\sup_{t\in[0,T]}
\biggl|
\int_W h(w) \mu\ssN_t(dw\times[y,1])-\int_W h(w) \mu_t(dw\times[y,1])
\biggr|^{2p}}
\\ \dsp {}
\le \frac{(1+2^{2p-1}C_h)C}{N^{2p\delta'}}+
\frac{2^{2p-1}C_h}{N^{2p}}\,,\ \ N\in\nintegers,
\arre
\end{equation}
where the constants in the right hand side is as in \thmu{SRPLLN}.
Since $2p\delta'>1$, we see that
\[
\EE{\sum_{N=1}^{\infty}\sup_{t\in[0,T]}
\biggl|
\int_W h(w) \mu\ssN_t(dw\times[y,1])-\int_W h(w) \mu_t(dw\times[y,1])
\biggr|^{2p}}<\infty,
\]
hence, in particular, we have
\begin{equation}
\eqna{HDLposdepprf4}
\limf{N}
\sup_{t\in[0,T]}
\biggl|
\int_W h(w) \mu\ssN_t(dw\times[y,1])-\int_W h(w) \mu_t(dw\times[y,1])
\biggr|=0,
\end{equation}
with probability $1$, which proves \eqnu{HDLposdepconv}.

Next we prove uniform almost sure convergence of $Y\ssN_i$ to $Y_i$
for $i=1,2,\ldots,L$.
As an analogy to 
\eqnu{latestarrivaltimeSItilde} and \eqnu{flowstarttime}, define
\begin{equation}
\eqna{jumptimelimit}
\arrb{l}\dsp
\tau_{i,0}=0,
\\ \dsp
\tau_{i,k+1}=\inf\{t>\tau_{i,k}\mid 
\nu_i(
\{(s,\xi)\in(\tau_{i,k},T]\times\preals\mid
\\ \dsp \phantom{\tau_{i,k+1}=}
0\le \xi\le w_i(Y_i(s-),s)\})>0\},\ k\in\pintegers,
\arre
\end{equation}
and
 \begin{equation}
\eqna{flowstarttimelimit}
\arrb{l}\dsp
\gamma_i(t)=\left\{\arrb{ll}\dsp 
(y_i,0),& 0\le t<\tau_{i,1},
\\ \dsp 
(0,\tau_{i,k}),& \tau_{i,k}\le t<\tau_{i,k+1},
\ k=1,2,\ldots.
\arre\right.
\arre
 \end{equation}
Comparing \eqnu{yc} and \eqnu{taggedlimit},
we have, with similar argument for 
\eqnu{YNthetaieqYNCfg},
\begin{equation}
\eqna{YieqyC}
Y_i(t)=y_C(\gamma_i(t),t),\ t\in[0,T].
\end{equation}

Quantities corresponding to
\eqnu{omegawedge}, \eqnu{omegavee},
\eqnu{MRUjumpseqdiffnu}, and \eqnu{difftime} are
\[ \arrb{l}\dsp
\breve{w}\ssN_{i,\wedge}(t)=
 w_i(Y\ssN_i(t-),t)\wedge w_i(Y_i(t-),t),t),
\\ \dsp
\breve{w}\ssN_{i,\vee}(t)=
 w_i(Y\ssN_i(t-),t)\vee w_i(Y_i(t-),t),t),
\\ \dsp
\tilde\mathcal{K}\ssN_i(0,t)=\{ \omega\in\Omega\mid 
\nu_i(\{(s,\xi)\mid
\breve{w}\ssN_{i,\wedge}(s)<\xi\le \breve{w}\ssN_{i,\vee}(s),
\ s\in(0,t]\ \}) =0  \},
\\ \dsp
\tilde{\sigma}\ssN_i=\inf\{t\in[0,T]\mid \tilde\mathcal{K}\ssN_i(0,t)^c\}.
\arre\]

A proof now proceeds with argument similar to that in \secu{eventcmp}.
An argument similar to that for \eqnu{recintdec} leads to
a bound
\begin{equation}
\eqna{Ydiffctbdlimit}
\arrb{l}\dsp
\tilde\mathcal{K}\ssN_i(0,t)^c
\\ \dsp {}
\subset
\{ \omega\in\Omega\mid \nu_i(\{(s,\xi)\mid
0\le \xi-\breve{w}\ssN_{i,\wedge}(s)
\le  C_W |Y\ssN_i(s)-Y_i(s)|,\ 
\\ \dsp \phantom{{}\subset\{ \omega\in\Omega\mid \nu_i(\{(s,\xi)\mid}
s\in(0,t]\ \})>0 \}.
\arre
\end{equation}
Since $\nu_i$ is a Poisson random measure, it holds 
with probability $1$ that $\nu_i(A)<\infty$
for a Borel set $A\subset [0,T]\times\preals$ of finite area.
Hence for almost all $\omega\in\Omega$
there exists $\eps=\eps(\omega)>0$ such that
\begin{equation}
\eqna{HDLposdepprf5}
\nu_i(\{(s,\xi)\mid
0\le \xi-\breve{w}\ssN_{i,\wedge}(s)\le \eps,\ s\in(0,t]\ \})=0.
\end{equation}
On the other hand,
applying the argument from \eqnu{HDLposdepprf3} to \eqnu{HDLposdepprf4},
\eqnu{YCNyCepsilon} implies
\begin{equation}
\eqna{HDLposdepprf10}
\limf{N} \sup_{(\gamma,t)\in\Delta_T}
|Y\ssN_C(\gamma,t)-y_C(\gamma,t)|=0,\ a.s..
\end{equation}
Therefore, for almost all $\omega\in\Omega$,
there exists an integer $N_0=N_0(\omega)$ such that
for $N\ge N_0$,
\begin{equation}
\eqna{HDLposdepprf6}
\sup_{(\gamma,t)\in\Delta_T}
|Y\ssN_C(\gamma,t)-y_C(\gamma,t)|\le \frac{\eps}{C_W}\,,\ N\ge N_0\,.
\end{equation}
Combining 
\eqnu{HDLposdepprf5}, \eqnu{HDLposdepprf6},
\eqnu{YNieqYNCfg}, \eqnu{YieqyC}, and
\eqnu{Ydiffctbdlimit},
\begin{equation}
\eqna{HDLposdepprf7}
\tilde\mathcal{K}\ssN_i(0,T)^c=\emptyset,\ N\ge N_0\,.
\end{equation}

Next, \eqnu{YNieqYNCfg} and \eqnu{YieqyC} imply
\[ \arrb{l}\dsp
|Y\ssN_i(t)-Y_i(t)|
\\ \dsp {} \le 
|Y\ssN_C(\gamma\ssN_i(t),t)-y_C(\gamma\ssN_i(t),t)|
+|y_C(\gamma\ssN_i(t),t)-y_C(\gamma_i(t),t)|
\\ \dsp {} \le 
\sup_{\gamma\in\Gamma_t}
|Y\ssN_C(\gamma,t)-y_C(\gamma,t)|
+|y_C(\gamma\ssN_i(t),t)-y_C(\gamma_i(t),t)|.
\arre \]
Comparing 
\eqnu{flowstarttimeSRP} and \eqnu{flowstarttimelimit}
we further have
\begin{equation}
\eqna{HDLposdepprf2}
\arrb{l}\dsp
|Y\ssN_i(t)-Y_i(t)|
\\ \dsp {} \le 
\sup_{\gamma\in\Gamma_t}
|Y\ssN_C(\gamma,t)-y_C(\gamma,t)|
\\ \dsp \phantom{\le} 
+|y_C((y\ssN_i,0),t)-y_C((y_i,0),t)|\chrfcn{t<\tau_{i,1}}
+\chrfcn{\tilde\mathcal{K}\ssN_i(0,t)^c}\,.
\arre
\end{equation}
Substituting \eqnu{HDLposdepprf7} in \eqnu{HDLposdepprf2},
\begin{equation}
\eqna{HDLposdepprf8}
\arrb{l}\dsp
|Y\ssN_i(t)-Y_i(t)|
\\ \dsp {} \le
\sup_{\gamma\in\Gamma_t}
|Y\ssN_C(\gamma,t)-y_C(\gamma,t)|
+|y_C((y\ssN_i,0),t)-y_C((y_i,0),t)|,
\ N\ge N_0\,.
\arre
\end{equation}
Since $y_C\in\Theta_T$,
$\limf{N} y\ssN_i=y_i$ implies
\begin{equation}
\eqna{HDLposdepprf9}
\limf{N} \sup_{t\in[0,T]} |y_C((y\ssN_i,0),t)-y_C((y_i,0),t)|=0,\ a.s..
\end{equation}
Combining 
\eqnu{HDLposdepprf10}, \eqnu{HDLposdepprf8}, \eqnu{HDLposdepprf9}
we have
\[
\limf{N} \sup_{t\in[0,T]} |Y\ssN_i(t)-Y_i(t)|=0,\ a.s..
\]
Therefore the almost sure uniform convergence of tagged
particle system holds,
which completes a proof of \thmu{HDLposdep}.

\appendix

\section{Point process with last-arrival-time dependent intensity.}
\seca{ptproc}

We will summarize the definition and basic formulas of the
point processes with last-arrival-time dependent intensities.
See \cite[\S3]{fluid14} and \cite[\S1.2]{fluid14k} for a proof.

In accordance with \cite{fluid14,fluid14k},
we will denote the point process with last-arrival-time dependent intensity
by $N=N(t)$, while we wrote $\tilde{\nu}$ in the main body of the
present paper to keep the symbol $N$ for the particle number.

Let $N=N(t)$, $t\ge0$, be a non-decreasing, right-continuous, 
non-negative integer valued stochastic process
on a measurable space with $N(0)=0$,
and for each non-negative integer $k$ define its $k$-th arrival time $\tau_k$
by
\begin{equation}
\eqna{arrivaltime}
\tau_k=\inf\{t\ge 0\mid N(t)\ge k\},\ \ k=1,2,\ldots,
\ \mbox{ and }\ \tau_0=0.
\end{equation}
The arrival times $\tau_k$ are non-decreasing in $k$, 
because $N$ is non-decreasing,
and since $N$ is also right-continuous, 
the arrival times are stopping times;
if we denote the associated filtration by 
$\dsp \mathcal{F}_t=\sigma[N(s),\ s\le t]$,
then $\{\tau_k\le t\} \in \mathcal{F}_t$, $t\ge 0$.

Let $\omega$ be a non-negative valued bounded continuous function 
of $(s,t)$ for $0\le s\le t$, and for $k=1,2,\ldots$ assume that
\begin{equation}
\eqna{intensity}
\begin{array}{l}\dsp
\prb{t< \tau_k\mid \mathcal{F}_{\tau_{k-1}}}
=\exp(-\int_{\tau_{k-1}}^t \omega(\tau_{k-1},u)\,du)
\ \mbox{ on }\ t\ge \tau_{k-1}\,.
\end{array}
\end{equation}
In particular, \eqnu{intensity} with $k=1$ implies
\begin{equation}
\eqna{intensity0}
\prb{N(t)=0}=\prb{\tau_1>t}=\exp(-\int_0^t \omega(0,u)\,du),\ \ t\ge0.
\end{equation}

If $\omega$ is independent of the first variable, then
\eqnu{intensity} implies that $N$ is the (inhomogeneous) Poisson process
with intensity function $\omega$.
We are considering a generalization of the Poisson process
such that the intensity function depends on the latest arrival time.

A construction of the point process with last-arrival-time dependent intensity
goes as follows.
Let $\omega:\ \preals^2\to\preals$ be a non-negative valued
bounded continuous function of $(s,t)$ for $0\le s\le t$,
for which we aim to construct a process satisfying \eqnu{intensity}.

Let $\nu$ be a Poisson random measure on $\preals^2$, 
with unit constant intensity 
\begin{equation}
\eqna{stdPRM}
\EE{\nu([a,b]\times[c,d])}=(b-a)(d-c)\,\ \ b>a>0,\ d>c>0,\ \ k\in\nintegers.
\end{equation}

Define a sequence of hitting times $\tau_k$, $k\in\pintegers$, inductively by
\begin{equation}
\eqna{arrivaltimedef}
\begin{array}{l}\dsp
\tau_0=0,\ \mbox{ and }\ 
\\ \dsp
\tau_k=\inf\{t\ge \tau_{k-1}\mid  \nu(\{(\xi,u)\in \preals^2\mid
\\ \dsp\phantom{\tau_k=\inf\{}
0\le \xi \le \omega(\tau_{k-1},u),\ \tau_{k-1}<u\le t\})>0\ \},
\\ \dsp
\ k=1,2,\ldots.
\end{array}
\end{equation}

Note that the definition is \textit{not} equivalent to the \textit{wrong}
formula such as
$\dsp
\tau_k=\inf\{t\ge 0\mid  \nu(\{(\xi,u)\in \preals^2\mid
0\le \xi \le \omega(\tau_{k-1},u),\ 0<u\le t\})\ge k\ \}.
$
We are recursively adding $1$ new arrival after the last arrival using
the renewed intensity $\omega(\tau_{k-1},\cdot)$ in \eqnu{arrivaltimedef}.

$\tau_k$ in \eqnu{intensity} is defined  by \eqnu{arrivaltimedef},
and the process $N(t)$
is defined by the reciprocal relation to \eqnu{arrivaltime}:
\begin{equation}
\eqna{Poissonklikeprocess}
N(t)=\max\{k\in\pintegers\mid \tau_k\le t\},\ t\ge0.
\end{equation}
Since $N(t)$ and $\tau_k$ are samplewise non-decreasing in $t$
and $k$, respectively, 
\eqnu{Poissonklikeprocess} and \eqnu{arrivaltime} are equivalent.
Also, \eqnu{intensity} follows from 
\eqnu{arrivaltimedef}.

$\{\tau_k\le t \}$ is in 
\begin{equation}
\eqna{Ft}
\mathcal{F}_t:= \sigma[\nu(A);
\ A\in\mathcal{B}(\preals^2),\ A\subset \preals\times[0,t],\ k\in\nintegers],
\end{equation}
and consequently $N$ is adapted to $\{\mathcal{F}_t\}$.

In coupling the stochastic ranking process with
the flow driven stochastic ranking process,
we will need a representation of $N$ by
the stochastic integration with respect to $\nu$ in \eqnu{stdPRM}
which is,
\begin{equation}
\eqna{lastarrivaltimedepSI}
N(t)=\int_{s\in(0,t]} \int_{\xi\in[0,\infty)}
\chrfcn{\xi\in[0,\omega(\tau^*(s-),s-))}\nu(d\xi\,ds),\ t\ge0,
\end{equation}
where $\tau^*:\ \Omega\times\preals\to\preals$ is 
the last arrival time up to time $t$:
\begin{equation}
\eqna{lastarrivaltime}
\tau^*(t)=\tau_{N(t)}=\inf\{ s\ge0\mid N(t)=N(s) \} \in[0,t],
\end{equation}
which satisfies a stochastic integration equation
\begin{equation}
\eqna{latestarrivaltimeSI}
\tau^*(t)=\int_{s\in(0,t]} \int_{\xi\in[0,\infty)}
(s-\tau^*(s-))\,\chrfcn{\xi\in[0,\omega(\tau^*(s-),s-))}\nu(d\xi\,ds),
\end{equation}
from which \eqnu{lastarrivaltimedepSI} follows.

For $t\ge t_0$ put
\begin{equation}
\eqna{numberdepPoissonprf11}
\Omega(t_0,t)=\int_{t_0}^t \omega(t_0,u)\,du.
\end{equation}
We have explicit formula
\begin{equation}
\eqna{numberdepPoissondecayposdep}
\arrb{l}\dsp
\prb{N(t)=N(s)}=\sum_{k\ge0} \prb{\tau_k\le s,\ t<\tau_{k+1}}
\\ \dsp\phantom{\prb{}}
=\sum_{k\ge0} \int_{0=:u_k< u_{k-1}<u_{k-2}<\cdots<u_1<u_0\le s}
\\ \dsp\phantom{\prb{}=\sum_{k\ge0} \int}\times
e^{-\sum_{i=0}^{k-1} \Omega(u_{i+1},u_i)-\Omega(u_0,t)}\,
\biggl(\prod_{i=0}^{k-1} \omega(u_{i+1},u_i)\,du_i\biggr),
\arre
\end{equation}
and for $\dsp \norm{\omega}=\sup_{0\le s\le t\le T} |\omega(s,t)|$,
\begin{equation}
\eqna{st_dep_diff}
\begin{array}{l}\dsp
0\le -\pderiv{}{t}\prb{N(t)=N(s)} 
\le \norm{\omega} (\prb{N(t)=N(s)}-\prb{N(t)=N(0)}) \le \norm{\omega},
\\ \dsp
0\le \pderiv{}{s}\prb{N(t)=N(s)} 
\\ \dsp\phantom{0\le}
=\sum_{k\ge0} \int_{0=:u_k< u_{k-1}<\cdots<u_1<u_0\le s} w(u_0,s)\,
\\ \dsp\phantom{0\le=\sum_{k\ge0} \int} \times
e^{-\sum_{i=0}^{k-1} \Omega(u_{i+1},u_i)-\Omega(u_0,t)}\,
(\prod_{i=0}^{k-1} \omega(u_{i+1},u_i)\,du_i)
\\ \dsp\phantom{\prb{}}
\le \norm{\omega} \prb{N(t)=N(s)}
\le \norm{\omega},
\\ \dsp
0\le s<t.
\end{array}
\end{equation}

\section*{Acknowledgements.}
The author would like to thank Professor Seiichiro Kusuoka, Professor Kazumasa Kuwada, and Professor Naoto Miyoshi, for their interest in the present work and for giving the author opportunities for seminar talks at Okayama University, University of Tokyo, and Tokyo Institute of Technology, respectively. The author in particular thanks Prof.\ Kusuoka for discussions.
This work is supported by JSPS KAKENHI Grant Number 26400146
from Japan Society for the Promotion of Science,
and by Keio Gijuku Academic Development Funds from Keio University.

\end{document}